\input amstex
\documentstyle{amsppt}
\pagewidth{6.5in}
%% STANDARD ABBREVIATIONS
\let\bs\backslash
\let\del\partial
\let\p\prime
\let\ol\overline
\let\ul\underline
\define\alg{\text{\rm alg}}

\define\NN{{\Bbb N}}
\define\HH{{\Bbb H}}
\define\CC{{\Bbb C}}
\define\RR{{\Bbb R}}
\define\QQ{{\Bbb Q}}
\define\PP{{\Bbb P}}
\define\ZZ{{\Bbb Z}}
\define\CH{\text{CH}}
\define\A{\text{A}}
\define\Ext{\text{Ext}}

\define\coker{\text{coker}}
\define\Image{\text{\rm Image}}
\define\End{\text{End}}
\define\MHS{\text{MHS}}
\define\MHM{\text{MHM}}
\define\DR{\text{DR}}
\define\cH{\Cal H}
\define\cO{\Cal O}
\define\cF{\Cal F}
\define\C{\Cal C}
\define\cS{\Cal S}
\define\X{\Cal X}
\define\MMot{\Cal{MM}}
\define\MF{\Cal{MF}}

\define\M{\Cal M}
\define\DD{\text{D}}
\define\Spec{\text{Spec}}
\define\Der{\text{Der}}

\define\tX{\widetilde{X}}
\define\qaq{\quad\text{ and }\quad}
\define\FB{F_B}
\define\GrFB{Gr_{\FB}}

%\define\XSk{X(S)/k}
\define\XSk{X_{S}/k}
%\define\VSk{V(S)/k}
\define\VSk{V_{S}/k}
\define\XSkd{X'_{S'}/k}
\define\XS{X/S}
\define\VS{V/S}
\define\XSS{\X_S/S}

\define\XKKk{X_K(K)/k}
\define\XKK{X_K/K}
\define\Xk{X/k}
\define\XSd{X'/S'}
\define\XC{X/\CC}
\define\SC{S/\CC}

\define\XCC{X_\CC}

\define\Rk{R/k}
\define\crXS{c^r_{\XSk}}
\define\csVS{c^r_{\VSk}}
\define\crXSDR{c^r_{\XSk,DR}}
\define\crXKKDR{c^r_{X_K(K)/k,DR}}
\define\crXSd{c^r_{\XSkd}}
\define\crXSDRd{c^r_{\XSkd,DR}}

\define\FkM{F_{\psi/k}}
\define\FkDR{F_{\DR/k}}

\define\rkDRkn{rk_{\DR/k}^{(\nu)}}
\define\rkDRkdn{rk_{\DR/k^{\p}}^{(\nu)}}
\define\rkDRn{rk_{\DR}^{(\nu)}}

\define\SDDnWX#1#2{\langle #1,#2 \rangle^{(\nu)}_{W,X}} 
\define\SDnWX{\langle\ ,\ \rangle^{(\nu)}_{W,X}} 
\define\SDnWY{\langle\ ,\ \rangle^{(\nu)}_{W,Y}}

\topmatter
\title  Algebraic Cycles and Mumford-Griffiths Invariants\endtitle
\leftheadtext{J. D. Lewis \& S. Saito}
\rightheadtext{Mumford-Griffiths Invariants}
\author James D. Lewis$^\dagger$ {\rm and} Shuji Saito \endauthor
%%\affil   \endaffil
\address  Department of Mathematical Sciences, University of Alberta,
Edmonton, Alberta, Canada
  \endaddress
\email lewisjd\@gpu.srv.ualberta.ca\endemail
\address  Department of Mathematical Sciences, University of Tokyo
Komaba, Meguro,  Tokyo 153, Japan
  \endaddress
\email sshuji\@msb.biglobe.ne.jp \endemail
%% \dedicatory   \enddedicatory
\date   May  2007
    \enddate
\thanks $\strut^\dagger$ Partially supported by a grant from the Natural
Sciences and Engineering Research Council of Canada.
\endthanks
\abstract
Let $X$ be a projective algebraic manifold and 
let $\CH^r(X)$ be the Chow group of algebraic cycles of 
codimension $r$ on $X$, modulo rational equivalence. 
Working with a candidate
Bloch-Beilinson filtration $\{F^{\nu}\}_{\nu\geq 0}$
on $\CH^r(X)\otimes \QQ$ due to the second author, 
we construct a space of arithmetic Hodge theoretic 
invariants $\nabla J^{r,\nu}(X)$ and corresponding
map $\phi_{X}^{r,\nu} : Gr_{F}^{\nu}\CH^r(X)\otimes \QQ
\to \nabla J^{r,\nu}(X)$, and determine conditions 
on $X$ for which the  kernel 
and image of $\phi_{X}^{r,\nu}$ are  ``uncountably large''.
\endabstract
\toc
\widestnumber\head{10b}
\head \S0. Introduction\endhead
\head \S1. Notation  \endhead
\head \S2. Good filtrations \endhead
\head \S3. Arithmetic de Rham cohomology  \endhead
\head \S4. Capturing nontrivial classes via  $\phi_{X}^{r,\nu}$\endhead
\head \S5. Cycle classes in higher extension groups \endhead
\head \S6. A map from cycle classes to de 
Rham/Mumford-Griffiths  invariants
\endhead
\head \S7. Detecting nonzero classes with trivial 
de Rham/Mumford-Griffiths invariant 
\endhead
\head \S8. A class of examples with trivial Mumford-Griffiths invariant \endhead
\specialhead {} Appendix: A conjectural overview\endspecialhead
\specialhead {} References\endspecialhead
\endtoc
\endtopmatter

\document
\head\S0. Introduction\endhead

Let $X$ be a projective smooth variety over $\CC$ and 
let $\CH^r(X)$ be the Chow group of algebraic cycles of codimension $r$ on $X$
modulo rational equivalence. A fundamental problem is to search for 
a reasonable set of invariants (e.g. Hodge theoretic) that provides
us with a good understanding of the structure of $\CH^r(X)$.  
The first significant step towards this problem was taken by 
Griffiths (1969) who 
defined  Abel-Jacobi maps
$$ 
\rho^r_X : \CH_{\hom}^r(X) \to J^r(X)
$$
where $\CH_{\hom}^r(X) = \ker\{ \CH^r(X) \to H^{2r}(X,\ZZ)\}$ is the subgroup of 
those cycle classes that are homologically equivalent to zero,
and where $J^r(X)$ is the $r$-th intermediate Jacobian of $X$ in the sense 
of  Griffiths \cite{Gri}. If the Griffiths Abel-Jacobi map were an isomorphism, 
then there would not be much
to explore in the world of algebraic cycles. That the map
is not surjective follows from the work of
Griffiths (op. cit.), using Hodge theory and monodromy
arguments; that the kernel is far from injective in general
is a consequence of Mumford's seminal work \cite{Mu}.

An important observation is the following natural isomorphism discovered 
by Carlson \cite{Ca},
$$ 
J^r(X) \simeq \Ext^{1}_{\MHS}(\ZZ, H^{2r-1}(X,\ZZ(r))),
$$
where $\MHS$ denotes the category of graded polarizable $\QQ$-mixed Hodge 
structures introduced by Deligne \cite{De1}. This implies that an element of 
$\CH_{\hom}^r(X)/\ker(\rho^r_X)$ 
is detected by an extension in $\MHS$. One may then have a naive expectation
that there may be a secondary cycle class map from $\ker(\rho^r_X)$ to 
higher extension groups $\Ext^p_{\MHS}$, which fails due to the fact 
that $\Ext^p_{\MHS}=0$ for $p\geq 2$ (\cite{Bei} (Cor. 1.10)).  
(Note that  $\Ext^p_{\MHS}=0$ for $p\geq 2$ also follows from the fact that 
$\Ext^1_{\MHS}(A,-)$ is  exact, using Carlson's 
explicit description of $\Ext^1_{\MHS}$ (\cite{Ca}. See also \cite{As1} 2.5.)
It was A. Beilinson who had an innovative idea to remedy the situation. 
He postulated the existence of a category
(called the category of mixed motives) whose higher extension
groups capture all elements of Chow groups.
A more precise formulation is the following conjecture.
In this paper we only consider Chow groups with rational coefficients:
$$
\CH^r(X;\QQ)=\CH^r(X)\otimes\QQ.
$$

\proclaim{Conjecture 0.1} For every projective smooth variety $X$ over $\CC$,
there exists a canonical and functorial filtration (called the Bloch-Beilinson filtration)
$$
\CH^r(X;\QQ)=F^0\CH^r(X;\QQ)\supset F^1\CH^r(X;\QQ)\supset 
F^2\CH^r(X;\QQ) \supset\cdots
$$
such that the following formula holds for each integer $\nu\geq 0$:
$$ 
F^\nu \CH^r(X;\QQ)/F^{\nu+1}\CH^r(X;\QQ) \simeq
\Ext^\nu_{\MMot}(1,h^{2r-\nu}(X)(r)).
$$
Here $\MMot$ denotes the (still conjectural) category of mixed motives over 
$\CC$ which contains as a full subcategory Grothendieck's category 
$\M$ of (pure) motives over $\CC$, $h^\ast(X)(r)\in \M$ denotes the 
cohomological object with Tate twist associated to $X$ and 
$1=h^0(\Spec(\CC))$.
\endproclaim

In \cite{Ja} and \cite{Sa2} it is proven that the Bloch-Beilinson filtration
is unique if it exists under the assumption of  
Grothendieck's standard conjectures.
Several candidates for the Bloch-Beilinson filtration 
have been proposed. From these
we adopt the filtration 
$$ 
F^\nu_{B} \CH^r(X;\QQ)\subset  \CH^r(X;\QQ) \quad (\nu\geq 0)$$
defined in \cite{Sa2}, Def.(1-3). 
We recall the definition of this filtration in \S 2.
\bigskip
The main results of this paper can be explained as
follows. We introduce the spaces
$\nabla J^{r,\nu}(X)$ of \it Mumford-Griffiths invariants \rm in \S3, which are
defined in terms of arithmetic de Rham cohomology.  
It is given by the cohomology of the complex 
$$
\Omega_{\CC/\ol{\QQ}}^{\nu-1}\otimes F^{r-\nu+1}H_{\DR}^{2r-\nu}(\XC) 
@>\nabla>>
\Omega_{\CC/\ol{\QQ}}^{\nu}\otimes F^{r-\nu}H^{2r-\nu}_{\DR}(\XC) @>\nabla>>
\Omega_{\CC/\ol{\QQ}}^{\nu+1}\otimes F^{r-\nu-1}H_{\DR}^{2r-\nu}(\XC),
$$
where 
$H_{\DR}^q(\XC)$ is the de Rham cohomology of $X/\CC$ with
the Hodge filtration $F^p H_{\DR}^q(\XC)$
and $\nabla$ is the arithmetic Gauss-Manin connection.
Then there is a cycle map (Proposition 3.7)
$$
\phi_{X}^{r,\nu}: Gr^{\nu}_{F_{B}}\CH^r(X;\QQ)
\to \nabla J^{r,\nu}(X).
$$
By ``forgetting'' the Hodge filtration, we also have a map
$$\phi_{X,DR}^{r,\nu}: Gr^{\nu}_{F_{B}}\CH^r(X;\QQ)
\to \nabla DR^{r,\nu}(X),$$
where 
$\nabla DR^{r,\nu}(X)$ is the coarser space of de Rham invariants
given by the cohomology of the complex
$$
\Omega_{\CC/\ol{\QQ}}^{\nu-1}\otimes H_{\DR}^{2r-\nu}(\XC) @>\nabla>>
\Omega_{\CC/\ol{\QQ}}^{\nu}\otimes H^{2r-\nu}_{\DR}(\XC) @>\nabla>>
\Omega_{\CC/\ol{\QQ}}^{\nu+1}\otimes H_{\DR}^{2r-\nu}(\XC).
$$
The first main result affirms that 
under various assumptions on $X$, the image of $\phi_{X,DR}^{r,\nu}$ is large,
where ``large'' has a similar meaning to that in the following theorem
of Mumford \cite{Mu}:

\proclaim{Theorem}
Let $X$ be a projective smooth surface over $\CC$.
Let $\A_0(X)\subset \CH_0(X;\QQ)$ denote the subgroup of the classes of 
zero-cycles of degree zero. 
Assume $H^0(X,\Omega^2_{X/\CC})\not=0$. Then $A_0(X)$ is 
infinite dimensional, viz., it is impossible to find
$Y_1,\dots,Y_N$, proper smooth connected curves with a morphism
$f: Y:=\coprod_{i=1}^N Y_i \to X$ such that
$f_* : \A_0(Y) \to \A_0(X)$ is surjective.
\endproclaim

We will give the following refinement of the above result:

\proclaim{Theorem 0.2}
Let $X$ be a projective smooth variety of dimension $d$ over $\CC$ and 
consider
$$\phi_{X,DR}^{d,\nu}: F_{B}^{\nu}\CH_0(X;\QQ)
\to \nabla DR^{d,\nu}(X).$$
Assume that there exists a dominant rational map
$\pi:\tX \dashrightarrow X$ such that the hard Lefschetz conjecture
$B(\tX)$ for $\tX$ holds (see Notation (v)).
Assume further that 
$H^0(X,\Omega^\nu_{X/\CC})\not=0$ for an integer $\nu\geq 2$. 
Then it is impossible to find
$\{f_i:Y_i\to X\}_{i\in \NN}$, a countable set of morphisms of proper smooth 
varieties over $\CC$ such that $\dim(Y_i)\leq \nu-1$ and
$$
F_B^\nu \CH_0(X;\QQ) \subset \Image\biggl\{
\underset{i\in \NN}\to{\bigoplus}\A_0(Y_i) @>{f_*}>>\; \A_0(X)\biggr\} +
\ker(\phi_{X,DR}^{d,\nu}),
$$ 
where $f_*$ is induced by $f_i$ for $i\in \NN$.
\endproclaim

We recall that for a smooth projective variety $W$, $B(W)$ is known to
hold if $W$ is obtained under successive operations of products and
hypersurface sections starting from curves, surfaces, Abelian varieties,
smooth complete intersections.

\medbreak

Our next result is an analogue of Theorem 0.2 for algebraic equivalence on
a hypersurface.  
By a general hypersurface $X$ of dimension $d$ and of degree $m$, we mean a
hypersurface corresponding to a point in a certain nonempty Zariski 
open subset of the universal family
of hypersurfaces of degree $m$ in $\PP^{d+1}$.

\proclaim{Theorem 0.3}
Let $X\subset \PP^{d+1}$ be a smooth hypersurface of degree $m\geq 3$
over $\CC$. Put
$$
k=\biggl[\frac{d+1}{m}\biggr]\qaq r=d-k,\; \nu=2r-d=d-2k
$$ 
and assume $k\geq 1,\; \nu\geq 2$, and the numerical condition:
$$
k(d+2-k) + 1 - {m+k\choose k} \geq 0.
$$ 
Consider
$$\phi_{X,DR}^{r,\nu}: F_{B}^{\nu}\CH^r(X;\QQ)
\to \nabla DR^{r,\nu}(X).$$

Assume that $X$ is general. Then it is impossible to find
$\{f_i:Y_i\to X\}_{i\in \NN}$, a countable set of morphisms of proper smooth 
varieties over $\CC$ such that $\dim(Y_i)\leq k+\nu-1$ and 
$$
A^r(X)\subset \Image\biggl\{
\underset{i\in \NN}\to{\bigoplus}\A_k(Y_i) @>{f_*}>>\; \A^r(X)\biggr\}
+ \ker(\phi_{X,DR}^{r,\nu}),
$$ 
where $f_*$ is induced by $f_i$ for $i\in \NN$.
Here $\A^r(X)\subset \CH^r_{\hom}(X;\QQ)$ denotes the subgroup of 
the cycle classes algebraically equivalent to zero. 
\endproclaim

\noindent
\underbar{Remarks 0.3.1}. (1) We will see in Proposition 2.7 that
$\CH^r_{\hom}(X;\QQ)=F_{B}^{1}\CH^r(X;\QQ)=F_{B}^{\nu}\CH^{r}(X;\QQ)$ 
in Theorem 0.3. 
\medbreak\noindent
(2) The hypersurfaces in question in Theorem 0.3,
albeit very interesting geometrically, are of a low degree.
In particular, since $k\geq 1$ it follows that
$H^{0}(X,\Omega^d_{X/\CC}) = 0$.
\bigskip

For the proof of the above results, we will introduce an integral invariant 
$rk(\Xi)$ for $\Xi\subset A^r(X)$, called the rank of $\Xi$ modulo 
$\ker(\phi_{X,DR}^{r,\nu})$, which measures the ``size'' of the image of
$\Xi\cap F_B^\nu\CH^r(X;\QQ)$ under $\phi_{X,DR}^{r,\nu}$.
The main technical result in this direction is stated in Theorem 4.2,
whose proof is based on an arithmetic version of Salberger's duality 
pairing, and the essential ideas rely on the constructions in \cite{Sa1}.
By working with the coarser de Rham invariants, we arrive at
similar statements for the Mumford-Griffiths invariants. 
\medbreak

With regard to $\ker(\phi_{X}^{r,\nu})$, 
if $X=X_o\times_{\ol{\QQ}}\CC$ with $X_o$ smooth projective over $\ol{\QQ}$, 
we introduce arithmetic Hodge theoretic invariants 
$H^{r-\nu,r}_{\ol{\QQ}}(X_o/\ol{\QQ})\subset H^r(X,\Omega^{r-\nu}_{X/\CC})$
which serve as an obstruction to countability of 
$\ker(\phi_{X}^{r,\nu})$. We will deduce the following:

\proclaim{Theorem 0.4} 
Let $X_o$ be a smooth projective variety over $\ol{\QQ}$
and $X=X_o\times_{\QQ}\CC$. Assume that the K\"unneth components of the
diagonal class of $X$ are algebraic.
If $\nu \geq 2$ and if $H_{\ol{\QQ}}^{r-\nu,r}(X_o/\ol{\QQ})\ne 0$, then
there are an uncountable number of  classes 
in the kernel of $\phi_{X}^{r,\nu}$.
\endproclaim

We will show how to construct a class of examples for which the assumptions 
of Theorem 0.4 are satisfied.
The examples arise from products involving $1$-motives over $\ol{\QQ}$,
in particular the product of smooth projective curves defined over $\ol{\QQ}$, 
and an Abelian variety defined over $\ol{\QQ}$. 
The precise statements of the main results appear in \S7 and \S8.
We believe that these theorems provide the strongest results to date
on the properties of $\phi_{X}^{r,\nu}$. For instance, Theorem 0.2
generalizes and strengthens Theorem 4.1 in \cite{MSa2}.
We note that it was Asakura (\cite{As2}) who first found an example of 
a nontrivial cycle with trivial Mumford-Griffiths invariant. 
\bigskip

To provide the reader with a better understanding of the results of 
this paper, we briefly discuss the technique of
taking a $\ol{\QQ}$-spread. This technique abounds in a number of works 
(e.g. \cite{Gr-Gr}, \cite{Le2}, \cite{MSa1-2}, \cite{As1}). 
The basic idea of the construction is the following. 
Given $X$, smooth projective 
over $\CC$, we can find a $\ol{\QQ}$-spread, namely a smooth affine
variety $S$ over $\ol{\QQ}$ with 
$$ f\;:\; \X_S \to S,$$
a projective smooth morphism of varieties over $\ol{\QQ}$, 
with a morphism $\eta:\Spec(\CC) \to S$ whose image is the generic point of
$S$, such that $\X_S\times_S \Spec(\CC)$, the base change via $\eta$, is 
isomorphic to $X$.
Similarly, if $\xi \in \CH^r(X;\QQ)$ is given,
then after making a base change if necessary,
$\xi$ has a lifting $\tilde{\xi} \in \CH^r(\X_S;\QQ)$. 
We then have the following variational version of $\phi^{r,\nu}_X$:
$$
\phi^{r,\nu}_{\XSS} \;:\; F_B^\nu\CH^r(\XSS;\QQ) \to 
\nabla J^{r,\nu}(\XSS).$$
We may retrieve $\phi^{r,\nu}_X$ from $\phi^{r,\nu}_{\XSS}$ via the
base change by $\eta$ (see \S 3 for the details).
\noindent

A key to the proof of Theorem 0.4 is the following construction: Put
$$
\Lambda^{r,\nu}(\XSS):= 
\hom_{\MHS}(\QQ(0),H^{\nu}(S/\CC,R^{2r-\nu} f_{\ast}\QQ(r))),
$$
$$
\Xi^{r,\nu}(\XSS) := \Ext^1_{\MHS}(\QQ(0),
H^{\nu-1}(S/\CC,R^{2r-\nu} f_{\ast}\QQ(r))),
$$
where $S/\CC=S\times_{\ol{\QQ}} \CC$, which we identify with
its underlying complex manifold, and 
$H^{\bullet}(S/\CC,R^{\bullet} f_{\ast}\QQ(r))$
denotes the cohomology of $S/\CC$ with coefficients in a local system,
which is endowed in a canonical way with a mixed Hodge structure by
the theory of mixed Hodge modules (\cite{MSa1}). (Alternatively,
the reader can consult \cite{A} for a different point of
view.) Then one defines the natural maps
$$
\lambda^{r,\nu}_{\XSS}\;:\; 
F_B^\nu\CH^r(\XSS;\QQ) \to \Lambda^{r,\nu}(\XSS),
$$
$$ 
\epsilon^{r,\nu}_{\XSS}\;:\;
\ker(\lambda^{r,\nu}_{\XSS}) \to \Xi^{r,\nu}(\XSS). 
$$
Broadly speaking, the construction is given in terms of an extension class
in the category of arithmetic mixed Hodge modules.
We then construct the following commutative diagram
under a suitable assumption on $\XSS$ (see Proposition 6.3)
$$
\matrix
\FB^\nu\CH^r(\XSS;\QQ) &@>{\phi^{r,\nu}_{\XSS}}>>& \nabla J^{r,\nu}(\XSS)\\ 
&&\\
\downarrow\rlap{$\lambda^{r,\nu}_{\XSS}$}&&
{\buildrel {}_{\ \cap}\over \downarrow}\\
 \Lambda^{r,\nu}(\XSS) &@>{\Phi^{r,\nu}_{\XSS}}>>& 
\nabla J^{r,\nu}(\XSS)\otimes_{\ol{\QQ}}\CC. \\
\endmatrix\tag{0.5}
$$
Now the cycle classes in Theorem 0.4 are constructed in 
$\ker(\lambda^{r,\nu}_{\XSS})$ and captured by $\epsilon^{r,\nu}_{\XSS}$.
\bigskip
The question now is whether the map $\Phi^{r,\nu}_{\XSS}$ in (0.5)
always exists for $S$ affine. What we can say is that if one
is willing to forgo the Hodge filtration by replacing the
space of Mumford-Griffiths invariants $\nabla J^{r,\nu}(\XSS)$ in (0.5)
by the coarser space of de Rham invariants 
$\nabla DR^{r,\nu}(\XSS)$ (see Definition 3.2), then 
for affine $S$, the diagram
corresponding to (0.5) does exist, and one
subsequently obtains
stronger results for cycles with trivial de Rham invariant
(Theorem 7.2). We conjecture that for affine $S$, the map
$\Lambda_{\alg}^{r,\nu}(\XSS) @>{\Phi^{r,\nu}_{\XSS}}>> 
\nabla J^{r,\nu}(\XSS)\otimes_{\ol{\QQ}}\CC$ always exists, where
$\Lambda_{\alg}^{r,\nu}(\XSS)=\text{\rm Image}(\lambda^{r,\nu}_{\XSS})$.
(For a more precise statement, see Conjecture 6.4.)
Roughly speaking, the discussion in the
Appendix says that this conjecture is a consequence
of the Hodge conjecture. If we put $S$ = Spec$(\CC)$
and let 
$$
\align 
\nabla J_{\alg}^{r,\nu}(X) &:= 
\text{\rm Image}\big\{Gr_{F_{B}}^{\nu}\CH^{r}(X;\QQ)
@>{\phi_{X}^{r,\nu}}>> \nabla J^{r,\nu}(X)\big\}\\
\nabla DR_{\alg}^{r,\nu}(X) &:= 
\text{\rm Image}\big\{Gr_{F_{B}}^{\nu}\CH^{r}(X;\QQ)
@>{\phi_{X,DR}^{r,\nu}}>> \nabla DR^{r,\nu}(X)\big\},
\endalign
$$
then a consequence of our aforementioned conjecture is that the
natural map 
$$
\nabla J_{\alg}^{r,\nu}(X)  \to \nabla DR_{\alg}^{r,\nu}(X).\tag{0.6}
$$
is an isomorphism. Put differently, and roughly
speaking, the Hodge conjecture implies that to
compute the Mumford-Griffiths invariant of an algebraic cycle,
it is sufficient to compute its de Rham invariant.\footnote{Since
the submission of this paper, M. Saito shared with us his work
in \cite{MSa4}, where now the conjectural assumptions, 
including Conjecture 6.4,  are
eliminated. Indeed as one of the referees pointed out, our Appendix 
seems already to have served as inspiration for \cite{MSa4}, and that
our discussion is still valuable.}
\bigskip
%The use of de Rham cohomology invariants to detect interesting
%algebraic cycles is not entirely new. For instance,
%Theorem 4 in \cite{N}, which is deduced from the
%Hodge theoretic Theorem 3 in \cite{N} together with
%a weight argument, is really in itself, a statement about
%de Rham cohomology. Yet it is that Theorem 4 that is used
%to detect nontrivial algebraic cycles in \cite{N} (Theorems 1 \& 2).
%\bigskip
We are very grateful to Matt Kerr
for meticulously reading parts of a preliminary version of this paper,
and for providing useful comments; and in particular for
sharing his ideas in \cite{Ke2}. We are also grateful to
V. Srinivas for pointing out his earlier work in \cite{Sr}. Indeed
the ideas presented in section 4 of this paper
share with \cite{Sr} a common methodology. Finally
we want to thank the referees for impressing upon us the need
to improve the presentation of our paper, by offering
their numerous constructive comments.

\newpage

\head\S1. Notation\endhead

\medbreak\noindent
(i) All fields in this paper are considered
as subfields of $\CC$. For a field $k$, let $\C_{k}$ be the category 
of smooth projective varieties over $k$. We simply
write $\C$ for $\C_{k}$ when $k=\CC$.
\medbreak\noindent
(ii) Unless otherwise indicated, $X \in \C_{k}$ 
will denote a smooth  projective
variety of dimension $d$.  (Periodically 
we will remind the reader that
$\dim X = d$.) Then $\CH^r(X)$ is the
Chow group of codimension $r$ algebraic cycles on $X$ modulo
rational equivalence, and $\CH^r(X;{\QQ}) =
\CH^r(X) \otimes {\QQ}$. Further, we put $A^{r}(X) :=
\CH^{r}_{\text{\rm alg}}(X;\QQ) \subset \CH^{r}(X;\QQ)$ to be the subgroup of
cycles algebraically equivalent to zero.
\medbreak\noindent
(iii) If $A\subset \RR$ is a subring, we denote by
$A(r)$,  $(2\pi\sqrt{-1})^{r}A$ the corresponding Tate
twist.
\medbreak\noindent
(iv) For $Y\in \C_k$, we write $Y_\CC := Y\times_{k}\CC$. 
It is sometimes more convenient to use the notation $Y/{\CC}$ for $Y_\CC$, 
which we also identify with its underlying complex manifold. 
\medbreak\noindent
(v) We fix a Weil cohomology theory for $\C_{k}$:
$$
\C_k \to Vec\; ; \; X\to H^\bullet (X),
$$
where $Vec$ denotes the category of finite dimensional vector spaces 
over a fixed field of characteristic zero.
Typical examples are given by 
$$
X\to H^\bullet(X/\CC,\QQ)
\quad\text{(singular\ cohomology)\  and }\quad 
X\to H_{\DR}^\bullet(X/k):=\HH^\bullet(X,\Omega^\bullet_{X/k}).
$$
There are well-known standard conjectures with Hg$(X)\Rightarrow
B(X) \Rightarrow C(X)$ (where Hg$(X)$ means the Hodge conjecture
for $X$).
Let $X\in \C_{k}$ with $d = \dim X$ and let $L_{X}\in H^2(X)$ be a hyperplane 
class. Let $i\leq d$ be an integer.
With regard to the hard Lefschetz isomorphism,
$$
L_{X}^{d-i} : H^{i}(X) @>\sim>> H^{2d-i}(X),
$$
the hard Lefschetz conjecture $B(X)$ asserts that the inverse 
$$
\big(L_{X}^{d-i}\big)^{-1} : H^{2d-i}(X) @>\sim>> H^{i}(X),
$$
is algebraic cycle induced. One consequence of $B(X)$ is a weaker conjecture 
$C(X)$ which asserts that the diagonal class 
$[\Delta_{X}]\in H^{2d}(X\times X)$, has algebraic K\"unneth components:
$$ [\Delta_{X}]=\sum_{i+j=2d} [\Delta_{X}(i,j)]\in 
\bigoplus_{i+j=2d} H^{i}(X)\otimes H^{j}(X)
\quad\text{ with }\Delta_{X}(i,j) \in \CH^{d}(X\times X;\QQ).$$ 
The reader can consult \cite{K1} for more details.
\medbreak\noindent
(vi) Let $f : X\to S$ be a smooth projective morphism of
quasiprojective varieties over a base field $k$. We define
$\CH^{r}(X/S;\QQ)$ to be $\CH^{r}(X;\QQ)$, and that
``$/S$'' will only affect how the filtration 
on $\CH^{r}(X/S;\QQ)$ in Definition 2.8 is defined.

\bigskip
\newpage  
\head\S2. Good filtrations \endhead

Fix a base field $k\subset \CC$. Our goal is to work with a
given filtration satisfying a number of good properties. 
More specifically, we require the following:

\proclaim{Definition 2.1} 
A filtration $F^\bullet \CH^\bullet$ of Bloch-Beilinson type on $\C_k$ is 
given by the following data:
For all $X\in \C_k$ and all $r$, there is a descending filtration
$$
\CH^r(X;\QQ) = F^0 \supset F^1\supset \cdots \supset 
F^{\nu}\supset F^{\nu +1}\supset\cdots
$$
which satisfies the following
\medskip
\noindent
{\rm (i)} $F^0 = \CH^r(X;\QQ)$ and $F^1 = \CH^r_{\hom}(X;\QQ)$.
\medskip
\noindent
{\rm (ii)} $F^{\nu}$ is preserved under the action of correspondences:
For $V,X\in \C_k$ and for $\Gamma\in \CH^q(V\times X;\QQ)$,
$$\Gamma_*(F^\nu\CH^s(V;\QQ))\subset F^\nu\CH^r(X;\QQ),$$
where
$\Gamma_* : \CH^s(V;\QQ) \to \CH^r(X;\QQ)$ with $s=r-q+\dim(V)$
is given by the formula
$$ 
\Gamma_\ast(\alpha)=(\pi_X)_\ast((\pi_V)^\ast(\alpha)\bullet\Gamma)
\quad \text{\rm for}\ \alpha\in \CH^s(V;\QQ),
$$
where $\pi_X:V\times X \to X$ and $\pi_V:V\times X \to V$ are the projections.
\medskip
\noindent
{\rm (iii)} The property {\rm (ii)} implies that we have the induced map
$$ Gr_F^\nu \Gamma_* : Gr_F^\nu\CH^s(V;\QQ) \to Gr_F^\nu\CH^r(X;\QQ).$$
Then $Gr_F^\nu \Gamma_*$ is the zero map if so is $\varphi_{\Gamma}^{2r-\nu}$ 
where 
$\varphi_{\Gamma}^{i} : H^{i-2(r-s)}(V) \to H^i(X)$ with $s=r-q+\dim(V)$
is given by the formula
$$ 
\varphi_{\Gamma}^{i}(\beta)=(\pi_X)_\ast((\pi_V)^\ast(\beta)\cup [\Gamma])
\ \text{\rm for}\ \beta\in H^{i-2(r-s)}(V),
$$
with $[\Gamma]\in H^{2q}(V\times X)$, the cohomology class of $\Gamma$.
\endproclaim
\medbreak

We point out the following consequences of Definition 2.1.

\proclaim{Proposition 2.2} Let $F^\bullet \CH^\bullet$ be a filtration 
of Bloch-Beilinson type on $\C_k$.
\medbreak\noindent
{\rm (i)} 
$F^\nu\CH^\bullet$ is preserved under push-forwards $f_*$ and pull-backs
$f^*$ for a morphism $f:X\to Y$ in $\C_k$.
\medbreak\noindent
{\rm (ii)} 
Let $X\in \C_k$ with $d=\dim(X)$.
Assume $C(X)$ holds and let $\Delta_X(p,q) \in \CH^{d}(X\times X;\QQ)$ be
as in \S1 (v). Then
$$
\Delta_X(2d-2r+\ell,2r-\ell)_{\ast}\biggl\vert_{Gr_F^{\nu}\CH^r
(X;\QQ)}
= \delta_{\ell,\nu}\cdot\text{\rm Identity}.
$$
\medbreak\noindent
{\rm (iii)} 
Under the same assumption as {\rm (ii)},
$$F^{\nu}\CH^r(X;\QQ)\bullet F^\mu\CH^s(X;\QQ) \subset 
F^{\nu +\mu}\CH^{r+s}(X;\QQ),$$
where $\bullet$ is the intersection product.
\endproclaim

\demo{Proof} Parts (i) and (ii) follows immediately from the definition. 
Part (iii) is \cite{Sa2}, Theorem(0-2).\qed
\enddemo

The following conjecture is due to Bloch and Beilinson:

\proclaim{Conjecture 2.3} There exists a filtration $F^\bullet \CH^\bullet$
of Bloch-Beilinson type satisfying:
$$
\bigcap_{\nu\geq 0}F^{\nu}\CH^r(X;\QQ)=0
\quad\text{ for all $X\in \C_k$ and all $r$}.
$$ 
\endproclaim
In general one can put 
$D^{r}(X;\QQ) := \cap_{\nu\geq 0}F^{\nu}\CH^r(X;\QQ)$ and
work modulo ``$D$-equivalence''.
We now introduce a specific filtration $F_B^\bullet \CH^\bullet$
of Bloch-Beilinson type, which is minimal among all
filtrations of Bloch-Beilinson type, namely it satisfies
$F_B^\bullet \CH^\bullet\subset F^\bullet \CH^\bullet$
for any filtration $F^\bullet \CH^\bullet$ of Bloch-Beilinson type. 
It is given in \cite{Sa2}, \S1.

\proclaim{Definition 2.4} 
For $\nu\geq 0$ we define $F_{B}^{\nu}\CH^r(X;\QQ)$ for all $X\in \C_k$
and for all $r\geq 0$ inductively as follows:
\bigskip
\item{\rm (1)} $F_{B}^{0}\CH^s(V;\QQ) = \CH^s(V;\QQ)$ for all $V\in \C_k$ 
and for all $s\geq 0$.
\bigskip
\item{\rm (2)} Assume that we have defined 
$F_{B}^{\nu}\CH^s(V;\QQ)$ for all $V\in \C_k$ and for all $s\geq 0$.
Then we define
$$
F_{B}^{\nu+1}\CH^{r}(X;\QQ)=\sum_{V,q,\Gamma}\text{\rm Image}
\big(\Gamma_\ast: F_{B}^\nu \CH^{r+d_V-q}(V;\QQ) \to \CH^r(X;\QQ)\big),
$$
where $V,q,\Gamma$ range over the following data:
\medskip
\item{\rm (a)} $V\in \C_k$ of dimension $d_V$,
\medskip
\item{\rm (b)} $r\leq q\leq r+d_V$,
\medskip
\item{\rm (c)} $\Gamma\in \CH^q(V\times X;\QQ)$ 
satisfying $\varphi_{\Gamma}^{2r-\nu}=0$, where
$\varphi_{\Gamma}^{2r-\nu}: H^{2s-\nu}(V) \to H^{2r-\nu}(X)$
is as in Def. 2.1 (iii).
\endproclaim
\medbreak

We have the following facts: 

\proclaim{Proposition 2.5} (\cite{Sa1} \S4) Let
$\rho^r_{\XCC} \;:\; \CH^r_{\hom}(X;\QQ)\to J^r(\XCC)\otimes\QQ$
be the Griffiths Abel-Jacobi map. Then
$$
\FB^2\CH^r(X;\QQ) \subset \ker(\rho^r_{\XCC})
\quad\text{ and }\quad
\FB^2\CH^r(X;\QQ)\cap \A^r(X)=\ker(\rho^r_{\XCC})\cap \A^r(X).
$$
\endproclaim
\medbreak

\proclaim{Proposition 2.6} (\cite{Sa2}, Thm. (1-1))
Assume $C(X)$ (cf. \S1 (v)) for all $X\in \C_k$ and that Conjecture 2.3 holds.
Then $F_B^\bullet\CH^\bullet$ is the only filtration of Bloch-Beilinson type.
\endproclaim
\medbreak

\proclaim{Proposition 2.7} 
Let $X\subset \PP^{N}$ be a smooth complete intersection of dimension $d$. 
For integers $r,\nu\geq 1$ we have:
$$
Gr_{F_B}^\nu \CH^r(X;\QQ)=0 \quad\text{ if }2r-\nu\not=d.
$$
\endproclaim
\demo{Proof} 
By Lefschetz theory we have for $0\leq i\leq 2d$:
$$
H^i(X)=
\left.\left\{
\aligned
 \QQ \cdot [Y]\\
  0
\endaligned\right.\qquad
\aligned
& i\not=d,\; \text{even}\\
& i\not=d,\; \text{odd}
\endaligned\right.
\leqno(*)
$$
where for $i=2m$, even, 
$j:Y\hookrightarrow X$ is the section by a general linear subspace
of codimension $m$ and $[Y]\in H^{2m}(X)$ denotes 
its cohomology class. Put $i=2r-\nu$ and assume $i\geq 1$ and $i\not=d$.
The diagonal $\Delta\subset X\times X$ induces
$$
\Delta_* : \GrFB^\nu\CH^r(X;\QQ) \to \GrFB^\nu\CH^r(X;\QQ)
\quad  \text{\rm and}\quad
\varphi_{\Delta}^{2r-\nu} : H^{2r-\nu}(X) \to H^{2r-\nu}(X)
$$
and both maps are the identity.
If $i$ is odd, then $\varphi_{\Delta}^{2r-\nu}=0$ by $(*)$ and hence
$\Delta_* =0$ by the definition of the filtration $\FB^{\bullet}$.
If $i=2m$, even, then 
$j_*: H^0(Y) @>\simeq>> H^{2m}(X)$ for $j:Y\hookrightarrow X$ as in $(*)$.
Since $B(X)$ holds for smooth complete intersections, we can find
$\Gamma\in \CH^{d_Y}(X\times Y;\QQ)$ ($d_Y=\dim(Y)=d-m$) such that
the induced map
$\varphi_{\Gamma}^{0} : H^{2m}(X) \to H^{0}(Y)$
is the inverse of the above map. Consider the induced maps
$$
\GrFB^\nu\CH^r(X;\QQ) @>{\Gamma_*}>> \GrFB^\nu\CH^{r-m}(Y;\QQ) @>{j_*}>> 
\GrFB^\nu\CH^r(X;\QQ),
$$
$$
H^{2m}(X) @>{\varphi_{\Gamma}^{0} }>> H^0(Y) @>{j_*}>> H^{2m}(X).
$$
The composite of the maps in the second row is the identity and so is 
that in the first row by the definition of the filtration $\FB^{\bullet}$.
Noting $r-m=\nu/2<\nu$, \cite{Sa2} Cor.(3-2) implies
$\GrFB^\nu\CH^{r-m}(Y;\QQ)=0$ which needs the fact that $B(Y)$ holds.
This proves the desired assertion. \qed
\enddemo
\medbreak

There is a variational (or relative) version of the aforementioned filtration.
The definition is completely analogous to that given in 
\cite{Sa3}, Definition (2-1).
As a base we fix a localization $S$ of a smooth quasiprojective variety $S$ 
over a field and let $\C_S$ be the category of $f:X\to S$, smooth 
projective morphisms.\footnote{We remind the reader that a smooth
morphism is a submersion at every point.}
For $f:X \to S$ in $\C_S$, let $\Omega^{\bullet}_{\XS}$ be the de Rham
complex of (Zariski) sheaves of relative differential forms of $X$ over $S$.
We define the de Rham cohomology sheaf:
$$
\cH^{i}_{\DR}(\XS) := \RR^{i}f_* (\Omega^{\bullet}_{\XS}),
$$
which is the Zarisiki sheaf of $\cO_S$-modules.
\medbreak
There is a direct generalization of Definition 2.4 to
the relative case.

\proclaim{Definition 2.8} 
For all $r,\nu\geq 0$ and for all $X\in \C_S$,
we define $F_{B}^{\nu}\CH^r(X/S;\QQ)\subset \CH^r(X;\QQ)$ 
in the following inductive way:
\bigskip
\item{\rm (1)} $F_{B}^{0}\CH^s(V/S;\QQ) = \CH^s(V;\QQ)$ for all $V\in \C_S$ 
and for all $s\geq 0$.
\bigskip
\item{\rm (2)} Assume that we have defined 
$F_{B}^{\nu}\CH^s(V/S;\QQ)$ for all $V\in \C_S$ and for all $s\geq 0$.
Then we define
$$
F_{B}^{\nu+1}\CH^{r}(X/S;\QQ)=\sum_{V,q,\Gamma}\text{\rm Image}
\big(\Gamma_\ast: F_{B}^\nu \CH^{r+d_V-q}(V/S;\QQ) \to \CH^r(X;\QQ)\big),
$$
where $V,q,\Gamma$ range over the following data:
\medskip
\item{\rm (a)} $V\in \C_S$ of relative dimension $d_V$,
\medskip
\item{\rm (b)} $r\leq q\leq r+d_V$,
\medskip
\item{\rm (c)} $\Gamma\in \CH^q(V\times X;\QQ)$ 
satisfying the condition $\varphi_{\Gamma,DR/S}^{2r-\nu}=0$, where
$$
\varphi_{\Gamma,DR/S}^{i}:
\cH_{\DR}^{i-2(r-s)}(V/S) \to \cH_{\DR}^{i}(X/S) \quad (s=r-q+d_V),
$$ 
is the homomorphism of $\cO_S$-modules, which is given by the same formula 
as in Definition 2.1 (iii), though  using
$[\Gamma]\in H^0(S,\cH_{\DR}^{2q}(V\times X/S))$, 
the cohomology class of $\Gamma$ in the de Rham cohomology.
\endproclaim
\medbreak

We note the following functoriality of the above filtration:
Let $\tau: T\to S$ be a morphism of schemes which are localizations of 
smooth quasiprojective varieties over a field. For $f:X\to S$ in $\C_S$ 
let $X_T=X\times_S T \in \C_T$ be the base change via $\tau$. 
Then, we have
$\tau^*(F_B^\nu \CH^r(X/S;\QQ))\subset F_B^\nu \CH^r(X_T/T;\QQ)$
under the pull-back 
$\tau^*:\CH^r(X;\QQ) \to \CH^r(X_T;\QQ)$.
If $\tau:T=\Spec(k) \to S$ is the generic point of $S$,
this gives the compatibility of absolute and relative versions of 
the filtrations on Chow groups in Definition 2.4 and 2.8.

\bigskip
\newpage

\head\S3. Arithmetic de Rham cohomology\endhead

Many of the ideas presented here in this section are inspired by the
lectures of M. Green \cite{Gr}. 
Fix a base field $k\subset \CC$. Let $\cS_k$ be the category of 
the affine integral schemes which are localizations of smooth
schemes over $k$. For example, for a finitely generated separable extension
of fields $K/k$, $\Spec(K)$ is an object of $\cS_k$.

\medbreak
\noindent
{\it 3.1. Arithmetic Gauss-Manin connection.} 
We recall the definition of the arithmetic Gauss-Manin connection. 
Take $S=\Spec(R)\in \cS_k$ and let $f:X \to S$ be smooth projective. 
Let $\Omega^{1}_{\XSk}$ (resp. $\Omega^{1}_{\XS}$) be the Zariski sheaf of 
relative differential forms of $X$ over $k$ (resp. $S$) and put 
$\Omega^{p}_{\XSk}=\overset{p}\to{\wedge}\Omega^{1}_{\XSk}$ 
(resp. $\Omega^{p}_{\XS}=\overset{p}\to{\wedge}\Omega^{1}_{\XS}$). 
Note that $\Omega^{\bullet}_{\XSk}$ is a complex under $d$.
We define the de Rham cohomology groups:
$$
H^{i}_{\DR}(\XSk) := {\HH}^{i}(X_{Zar},\Omega^{\bullet}_{\XSk})=
H^{i}(S_{Zar},\RR f_*\Omega^{\bullet}_{\XSk}),
$$
$$
H^{i}_{\DR}(\XS) := {\HH}^{i}(X_{Zar},\Omega^{\bullet}_{\XS}) =
H^{i}(S_{Zar},\RR f_*\Omega^{\bullet}_{\XS}).
$$
Put
$$
\text{\rm Filt}^{m}\Omega^{p}_{\XSk} := \text{\rm Im}\biggl(\Omega_{\Rk}^{m}
\otimes \Omega^{p-m}_{\XSk}\to \Omega^{p}_{\XSk}\biggr).
$$
Then:
$$
\text{\rm Gr}^{m}\Omega^{p}_{\XSk} \simeq \Omega^{m}_{\Rk}\otimes
\Omega^{p-m}_{\XS};
$$
Moreover:
$$
0\to{\text{\rm Filt}^{1}\Omega^{\bullet}_{\XSk}\over
\text{\rm 
Filt}^{2}\Omega^{\bullet}_{\XSk}} \to {\Omega^{\bullet}_{\XSk}
\over \text{\rm Filt}^{2}\Omega^{\bullet}_{\XSk}} \to {\Omega^{\bullet}_{\XSk}
\over \text{\rm Filt}^{1}\Omega^{\bullet}_{\XSk}}\to 0
$$
$$
 |\wr\hskip1in ||\hskip1in |\wr
$$
$$
0\to \Omega^{1}_{\Rk}\otimes \Omega^{\bullet}_{\XS}[-1]\to
{\Omega^{\bullet}_{\XSk}
\over \text{\rm Filt}^{2}\Omega^{\bullet}_{\XS}}\to
\Omega^{\bullet}_{\XS}\to 0
$$
Taking hypercohomology, we get a natural connecting map:
$$
\nabla := \nabla_{\XSk}: H^{i}_{\DR}(\XS) \to 
\Omega^{1}_{\Rk}\otimes H^{i}_{\DR}(\XS),
$$
called the arithmetic Gauss-Manin connection. By imposing Leibniz' rule, viz.,
$$
\nabla (\omega\otimes e) = d\omega \otimes e + (-1)^{m}
\omega\otimes \nabla e,
$$
one extends $\nabla$ to:
$$
\nabla  : \Omega^{m}_{\Rk}\otimes H^{i}_{\DR}(\XS) \to
 \Omega^{m+1}_{\Rk}\otimes H^{i}_{\DR}(\XS).
$$
It satisfies the following properties, which are consequences of 
the fact that $\nabla$ is identified with $d_1$ of the spectral sequences 
(3.1.1) and (3.1.2) below by \cite{KO}.
\medbreak\noindent
(i) $\nabla^{2} = 0$ (viz., flat connection),
\medbreak\noindent
(ii) (Griffiths transversality)
$$
\nabla\biggl(F^{p}H^{i}_{\DR}(\XS) \biggr)\subset \Omega^{1}_{\Rk}
\otimes F^{p-1}H^{i}_{\DR}(\XS),
$$
where
$F^{p}H_{\DR}^{\bullet}(\XS) := {\HH}^{\bullet}(\Omega^{\bullet\geq p}_{\XS})
\subset {\HH}^{\bullet}(\Omega^{\bullet}_{\XS})$
is the usual Hodge filtration.
\medbreak

Note that $\Omega^{\bullet}_{\XSk}$ is a filtered complex using
Filt$^{m}$. The corresponding spectral sequence is:
$$
E_{1}^{p,q} = \Omega^{p}_{\Rk}\otimes H^{q}_{\DR}(\XS)
\Rightarrow H^{p+q}_{\DR}(\XSk),
\tag{3.1.1}
$$
with $d_{1} = \nabla_{\XSk}$ (\cite{KO}). 
This is really the Leray spectral sequence, which by Deligne 
(\cite{De3}), is known to degenerate $E_{2}$. 
Likewise, the analogous Leray spectral sequence 
$$
E_{1}^{p,q} = \Omega^{p}_{\Rk}\otimes F^{r-p}H^{q}_{\DR}(\XS)\Rightarrow 
%F^{r}H_{\DR}^{p+q}(\XSk)):=
{\HH}^{2r}(\Omega^{\bullet\geq r}_{\XSk})
\tag{3.1.2}
$$
degenerates at $E_{2}$. We denote the corresponding Leray filtrations by:
$$
{F}_{L}^{\nu}{\HH}^{\bullet}(\Omega^{\bullet}_{\XSk})\subset
{\HH}^{\bullet}(\Omega^{\bullet}_{\XSk}),\quad
{F}_{L}^{\nu}{\HH}^{\bullet}(\Omega^{\bullet\geq r}_{\XSk})\subset
{\HH}^{\bullet}(\Omega^{\bullet\geq r}_{\XSk}).
\tag{3.1.3}
$$ 
\medbreak

\proclaim{Definition 3.2} 
{\rm (i)}
We put
$\nabla J^{r,\nu}(\XS)= 
Gr^{\nu}_{F_{L}}{\HH}^{2r}(\Omega^{\bullet\geq r}_{\XSk})$,
called the space of Mumford-Griffiths invariants of $\XS$. 
It is given by the cohomology of:
$$
\Omega_{\Rk}^{\nu-1}\otimes F^{r-\nu+1}H_{\DR}^{2r-\nu}(\XS) @>\nabla>>
\Omega_{\Rk}^{\nu}\otimes F^{r-\nu}H^{2r-\nu}_{\DR}(\XS) @>\nabla>>
\Omega_{\Rk}^{\nu+1}\otimes F^{r-\nu-1}H_{\DR}^{2r-\nu}(\XS).
$$
\medbreak\noindent
{\rm (ii)}
We put
$\nabla DR^{r,\nu}(\XS)= 
Gr^{\nu}_{F_{L}}{\HH}^{2r}(\Omega^{\bullet}_{\XSk})$,
called the space of de Rham invariants. It is given by the cohomology of:
$$
\Omega_{\Rk}^{\nu-1}\otimes 
H_{\DR}^{2r-\nu}(\XS) @>\nabla>>
\Omega_{\Rk}^{\nu}\otimes 
H^{2r-\nu}_{\DR}(\XS) @>\nabla>>
\Omega_{\Rk}^{\nu+1}\otimes 
H_{\DR}^{2r-\nu}(\XS).
$$
\medbreak\noindent
{\rm (iii)}
In case $S=\Spec(\CC)$ we simply write
$\nabla J^{r,\nu}(X)=\nabla J^{r,\nu}(\XS)$ and
$\nabla DR^{r,\nu}(X)=\nabla DR^{r,\nu}(\XS)$.
\endproclaim

Note that there is a natural map
$\nabla J^{r,\nu}(\XS) \to \nabla DR^{r,\nu}(\XS)$,
where one forgets the Hodge filtration.

\bigskip
\noindent
{\it 3.3. Arithmetic cycle class map.}
 Let $f:X\to S=\Spec(R)$ be as before. Let 
$$
{\Cal K}_{r,X}^{M}
= {\Cal O}_{X}^{\times}\otimes\cdots\otimes {\Cal 
O}_{X}^{\times}/\langle \tau_{1}\otimes\cdots
\otimes \tau_{i}\otimes\cdots\otimes \tau_{j}\otimes\cdots\otimes 
\tau_{r}\ |\ \tau_{i}+\tau_{j}=1,\ i\ne j\rangle,
$$
be the Milnor $K$-sheaf of $X$, and put 
$$
\ul{\Cal K}_{r,X}^{M}
=\ \text{\rm Image}\biggl({\Cal K}_{r,X}^{M} \to K^{M}_{r}({k}(X))\biggr).
$$
where $k(X)$ is the function field of $X$.
By the results of Gabber (or M\"uller-Stach, Elbaz-Vincent, see \cite{E-M}) 
$$
\CH^{r}(X) \simeq H^{r}_{\text{\rm \rm Zar}}(X,\ul{\Cal 
K}_{r,X}^{M}) = {\HH}^{r}(\ul{\Cal K}_{r,X}^{M} \to 0\to 
0\to \cdots).
$$
By torsion considerations (Suslin), the natural map 
${\Cal K}^{M}_{r,X} \to \Omega^{\bullet\geq r}_{\XSk}[r]$ given by
$$
\{f_{1},\ldots,f_{r}\} \mapsto \bigwedge_{j}d\log f_{j},\ f_{j}\in
{\Cal O}_{X}^{\times},
$$
factors through a morphism of complexes
$$
\big(\ul{\Cal K}_{r,X}^{M} \to 0\to 
0\to \cdots\big) \to \Omega^{\bullet\geq r}_{\XSk}[r],
$$
($\ul{\Cal K}_{r,X}^{M}$ and $\Omega^{r}_{\XSk}$
are in degree $0$)
and hence determines a morphism
$$
\crXS \;:\; \CH^{r}(X;\QQ) \to {\HH}^{2r}(\Omega^{\bullet\geq r}_{\XSk}),
\tag{3.3.1}
$$
called the arithmetic cycle class map. The same construction
also appears in \cite{Gr}, based on the discussion on p. 68
in \cite{E-P} and Th\'eor\`em 5 in \cite{So}.
(The reader can also consult the treatment by El-Zein \cite{EZ}.)
We will also use 
$$
\crXSDR \;:\; \CH^{r}(X;\QQ) \to {\HH}^{2r}(\Omega^{\bullet}_{\XSk}),
\tag{3.3.2}
$$
which is the composite of $\crXS$ and 
${\HH}^{2r}(\Omega^{\bullet\geq r}_{\XSk})\to
{\HH}^{2r}(\Omega^{\bullet}_{\XSk})$, the natural map.
These maps satisfy functoriality with respect to a base change
$\pi:S'\to S$ in $\cS_k$, namely we have the commutative diagrams 
($X'=X\times_S S'$):
$$
\matrix
\CH^{r}(X;\QQ) &@>{\crXS}>>& {\HH}^{2r}(\Omega^{\bullet\geq r}_{\XSk})\\
&\\
\downarrow\rlap{$\pi^*$} && \downarrow\rlap{$\pi^*$}\\
\CH^{r}(X';\QQ) &@>{\crXSd}>>& {\HH}^{2r}(\Omega^{\bullet\geq r}_{\XSkd}),
\endmatrix
\quad
\matrix
\CH^{r}(X;\QQ) &@>{\crXSDR}>>& {\HH}^{2r}(\Omega^{\bullet}_{\XSk})\\
&\\
\downarrow\rlap{$\pi^*$} && \downarrow\rlap{$\pi^*$}\\
\CH^{r}(X';\QQ) &@>{\crXSDRd}>>& {\HH}^{2r}(\Omega^{\bullet}_{\XSkd}),
\endmatrix
\tag{3.4}
$$

\proclaim{Definition 3.5} We introduce  filtrations
of Leray type on $\CH^{r}(\XS)$ as follows:
$$
\FkM^{\nu}\CH^{r}(\XS) = 
(\crXS)^{-1}\big(F^{\nu}_{L}{\HH}^{2r}(\Omega^{\bullet\geq r}_{\XSk})\big),
$$
$$
\FkDR^{\nu}\CH^{r}(\XS) = 
(\crXSDR)^{-1}\big(F^{\nu}_{L}{\HH}^{2r}(\Omega^{\bullet}_{\XSk})\big).
$$
\endproclaim
\medbreak

By definition we have
$\FkM^{\nu}\CH^{r}(\XS) \subset \FkDR^{\nu}\CH^{r}(\XS)$.
\medbreak 

\proclaim{Lemma 3.6}
Let $\pi:S'\to S$ be a morphism in $\cS_k$ and put $X'=X\times_S S'$.
We have
$$\pi^*(\FkM^{\nu}\CH^{r}(\XS))\subset \FkM^{\nu}\CH^{r}(\XSd),
\quad
\pi^*(\FkDR^{\nu}\CH^{r}(\XS))\subset \FkDR^{\nu}\CH^{r}(\XSd).$$
If $\pi$ is finite etale, then
$$\pi_*(\FkM^{\nu}\CH^{r}(\XSd))\subset \FkM^{\nu}\CH^{r}(\XS),
\quad
\pi_*(\FkDR^{\nu}\CH^{r}(\XSd))\subset \FkDR^{\nu}\CH^{r}(\XS).$$
\endproclaim
\demo{Proof}
This follows from the corresponding functoriality for 
(3.1.1) and (3.1.2).\qed
\enddemo
\medbreak

Now recall $\{\FB^{\nu}\CH^{r}(\XS;\QQ)\}_{\nu\geq 0}$ 
as introduced earlier in Definition 2.4 of \S2.
\medbreak

\proclaim{Proposition 3.7}
For all $\nu\geq 0$, 
$\FB^\nu\CH^{r}(\XS;\QQ) \subset \FkM^{\nu}\CH^{r}(\XS)$.
Hence there exists a natural map
$$
\phi_{\XS}^{r,\nu} : \FB^\nu \CH^r(\XS;\QQ) \to \nabla J^{r,\nu}(\XS)
$$
satisfying
$\phi_{\XS}^{r,\nu}( \FB^{\nu+1} \CH^r(\XS;\QQ))=0$.
It is functorial for  morphisms $Y\to X$ of projective smooth schemes
over $S$, also for base change $S'\to S$.
\endproclaim
\demo{Proof}
This is shown by the same argument as \cite{Sa3}, Prop.(2-1):
Let $\Gamma\in \CH^q(V\times X;\QQ)$ be as Definition 2.8 (c).
It induces a commutative diagram: 
$$
\matrix
\CH^{s}(V;\QQ) &@>{\csVS}>>& {\HH}^{2s}(\Omega^{\bullet\geq r}_{\VSk})\\
\downarrow\rlap{$\Gamma_*$}&&\downarrow\rlap{$\Gamma_*$}\\
\CH^{r}(X;\QQ) &@>{\crXS}>>& {\HH}^{2r}(\Omega^{\bullet\geq r}_{\XSk}) \\
\endmatrix
$$
where the vertical maps respect $F^\nu_B$ on Chow groups 
and $F^\nu_L$ introduced in (3.1.3). The induced map
$$
Gr^\nu_{F_L}\Gamma_* : Gr^\nu_{F_L}{\HH}^{2s}(\Omega^{\bullet\geq r}_{\VSk})\to
Gr^\nu_{F_L}{\HH}^{2r}(\Omega^{\bullet\geq r}_{\XSk})
$$
is then identified with the map induced by
$ \Gamma_*: H_{\DR}^{2s-\nu}(\VS) \to H_{\DR}^{2r-\nu}(\XS) $
via the identification of the above spaces with
the cohomology of the complex in Definition 3.2 (i).
The proposition follows by the induction on $\nu$.
\qed
\enddemo
\medbreak

We let
$$
\phi_{\XS,DR}^{r,\nu} : \FB^\nu \CH^r(\XS;\QQ) \to 
\nabla DR^{r,\nu}(\XS),\tag{3.7.1}
$$
denote the composite of $\phi_{\XS}^{r,\nu}$ and 
$\nabla J^{r,\nu}(\XS)\to \nabla DR^{r,\nu}(\XS)$,
the natural map.

In the case $S=\Spec(\CC)$ and $X\in \C$ we get the maps
$$\phi_{X}^{r,\nu}:=\phi_{X/\CC}^{r,\nu}
\; :\; \FB^\nu \CH^r(X;\QQ) \to \nabla J^{r,\nu}(X).
\tag{3.7.2}$$
$$\phi_{X,DR}^{r,\nu}:=\phi_{X/\CC}^{r,\nu}
\; :\; \FB^\nu \CH^r(X;\QQ) \to \nabla DR^{r,\nu}(X).
\tag{3.7.3}$$
\medbreak

\noindent
\underbar{Remark 3.8}.
If $X\in \C$ and $\nu=1$, $\phi_{X}^{r,\nu}$ is related to 
the Griffiths Abel-Jacobi map 
$\rho^r_X$ as follows. It is the Griffiths construction of infinitesimal
invariant of normal functions. We recall
$$ 
J^r(X) = H^{2r-1}(X_{an},\CC)/\big(F^{r}H^{2r-1}(X_{an},\CC) + 
H^{2r-1}(X_{an},\ZZ(r))\big).
$$
We have the comparison isomorphism
$H^q(X_{an},\CC)\simeq H_{\DR}^{q}(X/\CC)$ preserving the Hodge filtrations
and the arithmetic Gauss-Manin connection $\nabla$ annihilates the image of 
the subspace $H^q(X_{an},\QQ(r))\subset H^q(X_{an},\CC)$ for every $r$. 
Hence $\nabla$ induces
$$ 
\tau :J^r(X) \to \nabla J^{r,1}(X),
$$
and one can check that $\phi_{X}^{r,1} =\tau\circ\rho^r_X$.

\bigskip
\newpage

\head\S4. Capturing nontrivial classes via  ${\phi}^{r,\nu}_{X}$
\endhead
\medbreak

Let $\cF$ denote the set of the subfields
$\QQ\subset k\subset \CC$ whose transcendental degree over $\QQ$
is at most countable (bear in mind that $\CC$ has uncountable
transcendence degree over $\QQ$).
We want to prove infinite dimensionality statements about the image of 
${\phi}^{r,\nu}_{X}$ (from (3.7.2)) under certain conditions on $X$. 
In this section, we will actually prove infinite 
dimensionality results for the de Rham invariants, which 
immediately imply the same for the Mumford-Griffiths invariants. 

First, we need to introduce 
some integral invariants. We pick a base field $k\in \cF$.

\proclaim{Definition 4.1} 
Let $X \in \C$ be a projective smooth variety of dimension $d$ over $\CC$. 
Fix an integer $\nu> 1$, $\Xi$ a subset of $A^{r}(X)$ (cf. \S1 (ii)).  
We define $\rkDRkn(\Xi)$ to be the least of the
integers $\mu$ for which the following holds: 
There exists 
$\{f_i:Y_{oi}\times_k\CC \to X\}_{i\in \NN}$, a countable set of 
morphisms in $\C$ with $Y_{oi}\in \C_k$ 
such that $\dim Y_{oi} \leq d-r+\mu$ for all $i$ and 
$$
\Xi \subset \text{\rm Image}\biggl\{\underset{i\in \NN}\to{\bigoplus}
\CH_{d-r}(Y_i)@>f_{\ast}>> CH^{r}(X)\biggr\} \ + \FkDR^{\nu+1}\CH^{r}(X/\CC),
$$
where 
$f_*$ is induced by $f_{i}$ for $i\in \NN$, 
$\FkDR^{\nu+1}\CH^{r}(X/\CC)$ is as in Definition 3.5 with
$S=\Spec(\CC)$, and $Y_i=Y_{oi}\times_k\CC$. We also define
$$
\rkDRn(\Xi) = \text{\rm min}\{\rkDRkn(\Xi)\ \big|\; k\in \cF\}.
$$
Note that $\rkDRkn(\Xi) \geq \rkDRkdn(\Xi)$ if $k,k^{\p}\in \cF$ and 
$k\subset k^{\p}$, and that
$$ 0\leq \rkDRn(\Xi) \leq \rkDRkn(\Xi) \leq r.$$
\endproclaim

Our main result is the following:

\proclaim{Theorem 4.2} Let $W, X\in \C$ 
with $\dim X = d$ and $\Gamma\in \CH^{r}(W\times X)$ 
be given. Write
$$
rk^{(\nu)}_{\DR}(\Gamma) := 
rk^{(\nu)}_{\DR}\big(\Gamma_{\ast}(A_{0}(W))\big),
$$
where $\Gamma_{\ast} : A_{0}(W) \to A^{r}(X)$ is induced by 
$\Gamma$. Then we have the implication 
$$
rk^{(\nu)}_{\DR}(\Gamma) \leq \nu-1 \Rightarrow \text{\rm Image}
(\varphi^{\nu}_{\Gamma}) \subset F^{1}H^{\nu}_{\DR}(W/\CC),
$$
with 
$\varphi_{\Gamma}^{\nu} : H_{\DR}^{2(d-r)+\nu}(X/\CC)\to H^{\nu}_{\DR}(W/\CC)$ 
as in Def. 2.1 (iii), and where $F^{\bullet}$ is the Hodge 
filtration.
\endproclaim

The proof of this theorem proceeds in the same way as 
\cite{Sa1}, \S7 and \S8. It hinges on the following key lemma, 
which is an arithmetic counterpart of \cite{Sa1}, Theorem 7.1:

\proclaim{Lemma 4.3} 
Let $W,X\in \C_{k}$, $\Gamma \in \CH^{r}(W\times X)$ be given with 
$d=\dim X$. Introduce $m=d-r$ in Definition 4.1. 
Let $\eta$ be the generic point of $W$ and put 
$K= k(\eta)$ and $X_{K} = X\times_{k}K$. 
Assume given a  $k$-rational point $a\in W$ and put
$$
\xi :=\Gamma\bullet(\eta\times X)- \Gamma\bullet(a\times X) \in \CH^{r}(X_{K}).
$$
where $\bullet$ denotes the intersection product.
Suppose that there exists $\{f_i:Y_i\to X\}_{i\in \NN}$, a countable set of 
morphisms in $\C_k$ such that $\dim Y_i \leq m+\nu -1$ for all $i$ and  
$$
\xi \in \Image\big(\underset{i\in \NN}\to{\bigoplus} 
\CH_{m}(Y_{i K};\QQ)@>{f_\ast}>> \CH^r(X_{K},\QQ)\big) + 
\FkDR^{\nu+1}\CH^r(X_K/K),\quad (Y_{i K}=Y_i\times_{k} K)
$$
where 
$f_\ast$ is induced by $f_i$ for $i\in \NN$ and 
$\FkDR^{\nu+1}\CH^{r}(X_K/K)$ is defined as in Definition 3.5 for $S=\Spec(K)$.
Then we have
$$
\Image(\varphi_{\Gamma}^{\nu}) \subset F^{1}H^{\nu}_{\DR}(W/k),
$$
where $\varphi_{\Gamma}^{\nu} : H_{\DR}^{2(d-r)+\nu}(X/k) \to
H^{\nu}_{\DR}(W/k)$ is as in Def. 2.1 (iii).
\endproclaim

\demo{Proof} 
The idea of the proof originates from Salberger's duality
interpretation of Mumford's theorem (cf. \cite{Bl} \S1 Appendix): 
Writing for $K = k(\eta)$ and defining
$$
H_{\DR}^{p}(K/k) = \lim_{\buildrel \to\over {U\subset 
W}}H^{p}_{\DR}(U),
$$
with $U$ ranging over all nonempty Zariski open subsets of $W$, put
$$
\Phi_{\DR}^{\nu}(W) = \frac{H^{\nu}_{\DR}(W/k)}{\ker \big(
H^{\nu}_{\DR}(W/k) \to H^{\nu}_{\DR}(K/k)\big)}.
$$
The idea is to use a pairing defined for any $X\in \C_k$:
$$
\SDnWX : 
\CH^{r}(X_{K};\QQ)\otimes H_{\DR}^{2m+\nu}(X/k) \to \Phi_{\DR}^{\nu}(W),
$$
which satisfies the following properties:
\medbreak\noindent
(i) 
Let $\xi\in \CH^{r}(X_{K};\QQ)$ be as in Lemma 4.3. Then 
$$
\SDDnWX{\xi}{\beta}=0\text{  for any } \beta\in H_{\DR}^{2m+\nu}(X/k)
\;\Rightarrow \;
\Image(\varphi_{\Gamma}^{\nu}) \subset F^{1}H^{\nu}_{\DR}(W/k).
$$
\medbreak\noindent
(ii) 
$\SDnWX$ annihilates $\FkDR^{\nu+1}\CH^r(X_K/K)$ placed on the left-hand side.
\medbreak\noindent
(iii) 
If $f:Y\to X$ is a morphism in $\C_k$, there is a commutative diagram
$$
\matrix
\SDnWY \; : &
\CH^{r}(Y_{K}) &\otimes& H_{\DR}^{2m+\nu}(Y/k) &\to& \Phi_{\DR}^{\nu}(W)\\
& \downarrow\rlap{$f_*$} && \uparrow\rlap{$f^*$} && \| \\
\SDnWX \; : &
\CH^{r}(X_{K}) &\otimes& H_{\DR}^{2m+\nu}(X/k) &\to& \Phi_{\DR}^{\nu}(W)\\
\endmatrix
$$
\medbreak\noindent
(iv) Recall $m=d-r$. Then
$\SDnWY$ is trivial if $\dim(Y)\leq m+\nu-1$.
\medbreak
Lemma 4.3 follows immediately from this.
To construct the pairing, we note the decomposition
$$
\HH^{2r}(X_{K},\Omega^{\bullet}_{\XKKk}) = \bigoplus_{p+q=2r}
H^{p}_{\DR}(K/k)\otimes_{K}H^{q}_{\DR}(X_K/K) =
\bigoplus_{p+q=2r}
H^{p}_{\DR}(K/k)\otimes_{k}H^{q}_{\DR}(X/k),\tag{4.3.1}
$$
which is a consequence of the degeneration at $E_{1}$ of the Leray spectral 
sequence 
$$
E_{1}^{p,q} = \Omega_{K/k}^{p}\otimes_{K}H^{q}_{\DR}(\XKK) 
\Rightarrow \HH^{p+q}(X_K,\Omega_{\XKKk}^{\bullet})
$$
following from the product structure $W\times X$ over $k$, and 
de Rham base change (alternatively, one may apply the K\"unneth formula to 
$\Omega^{\bullet}_{\XKKk} \simeq \Omega^{\bullet}_{K/k}\otimes_{k} 
\Omega_{\Xk}^{\bullet}$,
where the differential $d$ on the left hand side is induced by 
$d_{K/k}\otimes 1 + (-1)^{\bullet}(1\otimes d_{\Xk})$ on the right hand side).
Thus from the product structure $W\times X$ over $k$, we
have the Leray filtration
$$
F_{L}^{\nu}\HH^{2r}(X_{K},\Omega^{\bullet}_{\XKKk}) = \bigoplus_{p\geq 
\nu,p+q=2r}H^{p}_{\DR}(K/k)\otimes_{k}H^{q}_{\DR}(\Xk).
$$
Then by composing $\crXKKDR$ for $X_K$ in 3.3.2 with the projection
$$
\HH^{2r}(X_{K},\Omega^{\bullet}_{\XKKk})  \to 
H^{\nu}_{\DR}(K/k)\otimes_{k} H^{2r-\nu}_{\DR}(X/k),
$$
and by using the natural pairing
$$H^{2(d-r)+\nu}_{\DR}(X/k)\otimes H^{2r-\nu}_{\DR}(X/k)\to
H^{2d}_{\DR}(X/k)\simeq k,$$
we arrive at $\SDnWX$. The above properties are now verified 
by the same argument as that in \cite{Sa1} \S7.\qed
\enddemo
\medbreak

\demo{Proof of Theorem 4.2}
Let $W,X,\Gamma$ be as in 4.2 and assume
$\rkDRn(\Gamma)\leq \nu-1$. Since the
assertion does not specify a ground field, we
may select $k\in \cF$ so that $X$ is defined over $k$
and the following further conditions hold:
\medbreak\noindent
$(a)$ 
There exists $X_o,W_o\in \C_k$ and $\Gamma\in \CH^r(W_o\times X_o,\QQ)$
which give rise to $W,X,\Gamma$ in 4.2 respectively by the base change via
$k\to \CC$.
\medbreak\noindent
$(b)$ 
$W_o$ has a $k$-rational point $a$.
\medbreak\noindent
$(c)$ 
There exists $\{f_i:Y_{oi} \to X_o\}_{i\in \NN}$, a countable set of 
morphisms in $\C_k$ such that $\dim Y_{oi} \leq m+\nu -1$ for all $i$ and 
$$
\Gamma_*(A_0(W))\subset \Xi \ {\buildrel \text{\rm def}\over =}\ 
 \Image\biggl\{\underset{i\in \NN}\to{\bigoplus}
\CH_{m}(Y_i;\QQ)@>f_{\ast}>> \CH^r(X,\QQ)\biggr\} + \FkDR^{\nu+1}\CH^r(X/\CC).
\quad (Y_i=Y_{oi}\times_k \CC)
$$
\medbreak

Let $\eta\in W_o$ be the generic point and put $K=k(\eta)$ and 
$$
\xi :=\Gamma\bullet(\eta\times X_o)- \Gamma\bullet(a\times X_o) \in 
\CH^{r}(X_{oK};\QQ).
$$
Theorem 4.2 follows from Lemma 4.3 if we show:
$$
\xi \in \Image\biggl\{\underset{i\in \NN}\to{\bigoplus}\CH_{m}(Y_{oiK};\QQ)
@>f_{K\ast}>> \CH^r(X_{oK},\QQ)\biggr\} + \FkDR^{\nu+1}\CH^r(X_{oK}/K).
\quad(Y_{oiK}=Y_{oi}\times_k K)
\tag{*}
$$
Take an embedding $\iota: K\hookrightarrow\CC$ of $k$-algebras, which is
possible by the assumption $k\in \cF$, and let
$\xi_\CC\in \CH^{r}(X;\QQ)$ be obtained from $\xi$ by the base change.
By definition $\xi_\CC\in \Gamma_*(A_0(W))$ and hence 
$\xi_\CC\in \Xi$ by $(c)$. 
To deduce $(*)$ from this, we consider the commutative diagram
$$
\matrix
\underset{i\in \NN}\to{\bigoplus}\CH_m(Y_{oiK};\QQ) &@>{f_{K*}}>>& 
\frac{\CH_m(X_{oK};\QQ)}{\FkDR^{\nu+1}\CH_m(X_{oK}/K)} \\
\downarrow\rlap{$\iota^*$}&&\downarrow\rlap{$\iota^*$} \\
\underset{i\in \NN}\to{\bigoplus}\CH_m(Y_i;\QQ) &@>{f_{\CC*}}>>& 
\frac{\CH_m(X;\QQ)}{\FkDR^{\nu+1}\CH_m(X/\CC)} \\
\endmatrix
$$
It suffices to show the injectivity of the induced map
$\coker(f_{K*})\to \coker(f_{\CC*})$.
This is shown by the same argument as \cite{Sa1}, Lemma (2.9),
which we recall for convenience of the readers.
Let $\cS_k$ be as in \S3.
For $S\in \cS_k$ we put
$$ \Lambda_S=\coker\biggl\{ 
\underset{i\in \NN}\to{\bigoplus}\CH_m(Y_{oiS};\QQ) @>{f_{S*}}>>
\frac{\CH_m(X_{oS};\QQ)}{\FkDR^{\nu+1}\CH_m(X_{oS}/S)} \biggr\}.
$$
It suffices to show the injectivity of the induced map
$g^* : \Lambda_K \to \Lambda_S$ for any morphism $g: S\to \Spec(K)$ in $\cS_k$,
where for notational convenience $\Lambda_K$ means $\Lambda_{\Spec(K)}$.
We may assume $g$ is of finite type, for 
if $S=\Spec(L)$ for a finite field extension $L$ of $K$, then 
we have the functorial map
$g_* : \Lambda_L \to \Lambda_K$ such that $g_*g^*$ is  multiplication  
by the degree of $L/K$, which implies the desired injectivity.
Here we used Lemma 3.6. 
In general, by what is just shown, we may assume that there exists
a section $i: \Spec(K)\to S$ of $g$. Then the composite
$ \Lambda_K @>{g^*}>> \Lambda_S  @>{i^*}>> \Lambda_K$
is the identity and hence $g^*$ is injective.
This completes the proof of Theorem 4.2.\qed
\enddemo
\bigskip
%Note: One could ask about uncountable sums.
%E.g. $X/\CC$ a surface, $r=2$ and $p\in X(\CC)$
%varies. In this case $k=\CC$, but there
%is no embedding $k(X) = \CC(X) \hookrightarrow \CC$
%over $\CC$

Recall $\{\FB^{\nu}\CH^{r}(X;\QQ)\}_{\nu\geq 0}$ 
as introduced earlier in Definition 2.4 of \S2.
Theorem 0.2 follows from the following:

\proclaim{Corollary 4.4} Let  $X\in \C$ with $\dim X = d$. Assume that
either $\nu \leq 2$ or that there exists a dominant rational morphism 
$\tilde{X} \to X$ such that the condition $B(\tilde{X})$ holds (cf. \S1 (v)). 
Then 
$$ \rkDRn(\FB^{\nu}\CH_0(X;\QQ)) \leq \nu-1 \Rightarrow
H^{0}(X,\Omega_{X/\CC}^{\nu})=0.$$
\endproclaim
\demo{Proof} 
This is an immediate consequence of Theorem 4.2 applied to the case 
$W=X$ and $\Gamma=\Delta_X$, together with the following result
(\cite{Sa1}, Theorem (6.2)).\qed
\enddemo
\medbreak
\proclaim{Theorem 4.5} Let the assumption be as in 4.4.
Then there exists a morphism $f : Y\to X$ in $\C$ such that
$\dim Y \leq \nu-1 $ and 
$$ A_{0}(X) \subset \FB^{\nu}\CH_0(X;\QQ) + f_*(A_0(Y)).$$
\endproclaim
\medbreak

Theorem 0.3 follows from the following result:

\proclaim{Corollary 4.6}
Let $X\subset \PP^{d+1}$ be a general hypersurface of degree $m\geq 3$
over $\CC$, where ``general" means corresponding to a point in a certain 
nonempty Zariski open subset of the universal family of hypersurfaces of 
degree $m$ in $\PP^{d+1}$. Let 
$k = \big[\frac{d+1}{m}\big]$ (greatest integer) and 
$r=d-k$ and $\nu = d-2k$,
and assume that the following numerical condition holds:
$$
k(d+2-k) + 1 - {m+k\choose k} \geq 0.
$$
Then $\rkDRn(A^r(X))=\nu$
(Note that $A^r(X)\subset \FB^\nu\CH^r(X)$ by Proposition 2.7).
\endproclaim
\demo{Proof}
Let ${\Cal G}_{X}(k)$ be the moduli space of the linear subspaces of dimension 
$k$ contained in $X$. Then there exists a smooth projective subvariety 
$W \subset {\Cal G}_{X}(k)$, of dimension $\nu$, such that 
the induced cylinder map
$$
[\Gamma]_{\ast} : H^{d-2k}(W,\CC)\to H^{d}(X,\CC),
$$
is surjective, where $\Gamma \in \CH^{r}(W\times X)$ is
the incidence correspondence (\cite{Le3}). The numerical condition implies
$H^{k}(X,\Omega^{d-k}_{X/\CC})\ne 0$ (see \cite{Le1}) so that
$\Image([\Gamma]_{\ast}) \not\subset F^{k+1}H^{d}(X,\CC)$. 
It implies 
$\Image([\Gamma]^{\ast})\not\subset F^1H^{\nu}(W,\CC)$,
where 
$
[\Gamma]^{\ast} : H^{d}(X,\CC)\to H^{\nu}(W,\CC)
$ 
is the dual of $[\Gamma]_{\ast}$ (\cite{Sa1}, Cor.(8.4)).
It follows from Theorem 4.2 that $\rkDRn(\Gamma) = \nu$. 
In fact, it is also the case that 
$\Gamma_{\ast}\big(A_{0}(W)\big) = A^r(X)$ (\cite{Le3}), thus 
$\rkDRn\big(A^r(X)\big) = \nu$.\qed
\enddemo
\bigskip
\noindent
\underbar{Remark 4.7}. The assumption of $X$ being
general in Corollary 4.6 is needed in \cite{Le3}
to arrive at results on the level of Chow groups, in particular
our above assertion $\Gamma_{\ast}\big(A_{0}(W)\big) = A^r(X)$.
For instance, we need (among other things),
the nonsingularity and predicted dimension
of the Fano variety of $k$-planes on $X$.

\newpage

\head\S5. Cycle classes in higher extension groups
\endhead

In this section we give a brief explanation of cycle classes in higher 
extension groups in suitable categories following the works of 
Green-Griffiths, M. Asakura and M. Saito 
(cf. \cite{Gr}, \cite{As1} and \cite{MSa2}).
It is motivated by  Beilinson's conjectural formula in \S0. 
Fix a base field $k\subset \CC$ as in the previous sections.

\medbreak
\noindent
{\it 5.1. Mumford-Griffiths invariants as higher extensions.} 
Let $R$ be a localization of a smooth algebra over $k$ and let 
$\DD_{R}\subset \End_{k}(R)$ be the ring of differential operators 
generated as a subring by $\Theta_R=\Der_{k}(R,R)$ (derivations) and 
$R$ (scalars). If $R$ has a local coordinate $\{x_i\}_{1\leq i\leq n}$ and 
$\{\del_i\}_{1\leq i\leq n}$ is the dual basis of 
$\{d x_i\}_{1\leq i\leq n}$, we have
$$ 
\DD_{R} = \bigoplus_{\alpha\in \ZZ_{+}^{n}} \; R \cdot\del^\alpha
\quad (\del^\alpha:=\del_1^{\alpha_1}\del_2^{\alpha_2}\cdots
\del_n^{\alpha_n}).
$$
It is endowed with the filtration of differential order:
$$ 
F_p \DD_{R}= \bigoplus_{|\alpha|\leq p} \; R \cdot \del^\alpha
\quad (|\alpha|=\sum_{i=1}^n \alpha_i).
$$
A filtered $\DD_{R}$-module is a pair $(M,F)$ consisting of a $\DD_{R}$-module $M$ and 
an increasing filtration of {\it finite} $R$-modules $F_p M\subset M$ 
($p\in \ZZ$) satisfying
\bigskip
\item{(1)} $M=\bigcup_{p\in \ZZ} F_p M \quad (F_p M=0\quad p<<0),$
\bigskip
\item{(2)} $F_p \DD_{R} \cdot F_q M \subset F_{p+q} M$.
%\item{(3)} Repacing $\subset$ by $=$ in (2) for large $p,q$.
\bigskip
Let $\MF(R)$ be the category of filtered $\DD_{R}$-modules.
For an object $M=(M,F)$ of $\MF(R)$ its Tate twist $M(r)$ is defined to be
$(M,F(r))$ with $F(r)_p=F_{p-r}$. 
Let $R(m)\in \MF(R)$ be the $\DD_{R}$-module $R$ endowed with the filtration 
$F_p R =0$ if $p<m$ and $F_p R =R$ if $p\geq m$. 
$\MF(R)$ becomes an exact category by defining a complex in $\MF(R)$
$$ 
(M_1,F_1) \to (M_2,F_2) \to (M_2,F_2) 
$$
to be exact if and only if 
$$ 
Gr^{F_1} M_1 \to  Gr^{F_2} M_2 \to  Gr^{F_3} M_3 
$$
is an exact sequence of $R$-modules. 
Thus we can consider higher extension groups in the sense of Yoneda in 
$\MF(R)$.  

Let $f:X\to S=\Spec(R)$ be a smooth projective morphism. 
The de Rham cohomology $M= H^{q}_{\DR}(\XS) =
\HH^q(\Omega^{\bullet}_{\XS})$ 
with the filtration
$F_p M:= F^{-p} H^{q}_{\DR}(\XS)$ gives rise to an object of $\MF(R)$: We let 
$\theta\in \Theta_R$ act on $M$ via
$\nabla_{\theta}$, the covariant derivative of $\theta$ with respect to
the arithmetic Gauss-Manin connection
$$
\nabla : H^{q}_{\DR}(\XS)\to \Omega^{1}_{R/k}\otimes H^{q}_{\DR}(\XS).
$$

\proclaim{Proposition 5.2}  (cf. \cite{Gr} and \cite{As1}, 3.2) 
For integers $p,q\geq 0$, 
$\Ext_{\MF(R)}^p(R(0),H^q_{DR}(\XS)(r))$ 
is isomorphic to the cohomology of the following complex 
$$ 
\Omega^{p-1}_{R/k}\otimes F^{r-p+1}H^{q}_{\DR}(\XS)
@>\nabla>> \Omega^{p}_{R/k}\otimes F^{r-p}H^{q}_{\DR}(\XS)
@>\nabla>> \Omega^{p+1}_{R/k}\otimes F^{r-p-1}H^{q}_{\DR}(\XS).
$$
\endproclaim

\proclaim{Corollary 5.3} We recall that $S =  \Spec(R)$. 
For integers $r,\nu\geq 0$ there is 
a canonical isomorphism
$$ 
\nabla J^{r,\nu}(\XS)\simeq
\Ext_{\MF(R)}^\nu(R(0),H_{\DR}^{2r-\nu}(\XS)(r)).
$$
\endproclaim

\medbreak
\noindent
{\it 5.4. Arithmetic mixed Hodge modules.} 
In the above construction, the $\QQ$-structure is not taken into account. 
The theory of arithmetic mixed Hodge modules remedies the defect and 
gives a refinement of the Mumford-Griffiths invariants.

\medbreak
For a smooth variety $X$ over $k$, the category $\MHM(X)$ of mixed Hodge
modules on $X$ was defined by Morihiko Saito (\cite{MSa3}). 
Let $X(\CC)$ be the set of $\CC$-valued points over $k$ endowed with 
the usual analytic topology. An object of $\MHM(X)$, called a mixed Hodge 
module on $X$, is given by a triple
$(K_\QQ,(M,F),\iota)$ satisfying certain properties,
where $K_\QQ$ is a perverse sheaf of $\QQ$-vector spaces of finite rank 
on $X(\CC)$, and $(M,F)$ is a filtered
$D_X$-module which is holonomic and has regular singularities, and 
$\iota$ is a quasi-isomorphism (Riemann-Hilbert correspondence)
$$
DR(M) \simeq K_\QQ\otimes\CC\quad\text{ with }
DR(M)=M\otimes_k \Omega^{\bullet}_{X(\CC)} [\dim(X)].
$$
(To be more precise, a mixed Hodge module is equipped with a weight filtration.
Since it is not used in this paper, we have omitted it for the sake of 
simplicity. The reader may consult the brief exposition in \cite{AS1} \S2.)
The category $\MHM(X)$ of mixed Hodge
modules on $X$ is an abelian category. For a morphism $f:X\to Y$ there are the standard 
operations 
$f_*$, $f^*$, $f_!$, $f^!$, etc., on the derived category of 
bounded complexes of mixed Hodge modules on $X$ and $Y$. 
There exists a cycle class map
$$ \rho^r_X \;:\; \CH^r(X;\QQ) \to \Ext^{2r}_{\MHM(X)}(\QQ_X(0),\QQ_X(r)),$$
where $\QQ_X(r)$ is the Tate object in $\MHM(X)$ whose underlying perverse 
sheaf is $\QQ_X[\dim(X)]$ (the constant local system with a degree shift)
and whose underlying $D_X$-module is $\cO_X$ with the filtration given by
$F_{-r}\cO_X=\cO_X$ and $F_{-r-1}\cO_X=0$.

Let $f:X\to S$ be a morphism of quasiprojective smooth varieties over $k$.
The standard operations provide us with the Leray spectral sequence
$$ E_2^{p,q}= \Ext^p_{\MHM(S)}(\QQ_S(0),R^q f_*\QQ_X(r))
\Rightarrow  \Ext^{p+q}_{\MHM(X)}(\QQ_X(0),\QQ_X(r)),$$
where $R^qf_* : \MHM(X)\to \MHM(S)$ is the derived functor on the categories
of mixed Hodge modules. The decomposition theorem for mixed Hodge modules
implies that the spectral sequence degenerates at the $E_2$-term, provided $f$ is 
proper. 

In what follows we assume further that $f$ is proper and smooth. Letting
$$
F^\nu_L \Ext^{2r}_{\MHM(X)}(\QQ_X(0),\QQ_X(r))\subset
\Ext^{2r}_{\MHM(X)}(\QQ_X(0),\QQ_X(r))
$$
be the filtration inducing the spectral sequence, we can show
$$
\rho^r_X(\FB^\nu\CH^r(\XS;\QQ))\subset 
F^\nu_L \Ext^{2r}_{\MHM(X)}(\QQ_X(0),\QQ_X(r))
$$
by the same argument as in \cite{Sa3}, Prop.(2-1) 
(also see the proof of Prop.3.7). Hence we get the induced map
$$
\rho^{r,\nu}_{\XS} \;:\; Gr^\nu_F \CH^r(\XS;\QQ) \to 
\Ext^\nu_{\MHM(S)}(\QQ_S(0),R^{2r-\nu} f_*\QQ_X(r)).
$$

If $S=\Spec(R)$, there is a natural functor
$$
r_{\MF}\; : \; \MHM(S) \to \MF(R)
$$
which sends a mixed Hodge module to its underlying filtered $D_S$-module.
We have $r_{\MF}(\QQ_S(m))=R(m)$ and 
$r_{\MF}(R^qf_*\QQ_X(r))=H_{\DR}^q(\XS)(r)$
and the map $\phi^{r,\nu}_{\XS}$ in Proposition 3.7
coincides with the composite of $\rho^{r,\nu}_{\XS}$ and the map
induced by $r_{\MF}$:
$$
\Ext^\nu_{\MHM(S)}(\QQ_S(0),R^{2r-\nu} f_*\QQ_X(r)) \to
\Ext_{\MF(R)}^\nu(R(0),H_{\DR}^{2r-\nu}(\XS)(r))\simeq \nabla J^{r,\nu}(X/S),
$$
where the last isomorphism is given in Corollary 5.3.

Let $\eta:\Spec(\CC) \to S$ be a $k$-morphism and 
$X_\eta=X\times_S \Spec(\CC)$ be the base change.
Noting that $\MHM(\Spec(\CC))$ is the category $\MHS$ of graded polarized 
$\QQ$-mixed Hodge strucutres, $\eta$ induces a natural functor
$$
r_{H,\eta}\; : \; \MHM(S) \to \MHS.
\tag{5.4.1}
$$
%where $\MHS$ denotes the category of mixed Hodge structures.
It satisfies 
$$
r_{H,\eta}(\QQ_S(m))=\QQ(m)\ \text{\rm and}\ 
r_{H,\eta}(R^qf_*\QQ_X(r))=H^q(X_{\eta,an},\QQ(r)),
$$
and the composite of $\rho_X^{r,1}$ with the map
$$ 
\Ext^1_{\MHM(S)}(\QQ_S(0),R^{2r-1} f_*\QQ_X(r))  \to 
\Ext^1_{\MHS}(\QQ(0),H^{2r-1}(X_{\eta,an},\QQ(r)))\simeq J^r(X_\eta),
$$
coincides with the Griffiths Abel-Jacobi map for $X_\eta$, where
the first map is induced by $r_{H,\eta}$ and 
the second map is Carlson's isomorphism (\cite{Ca}).
\bigskip
\noindent
\underbar{Remark 5.5}. As pointed out by one of the
referees, the map induced by $r_{\MF}$ into
$\nabla J^{r,1}(X/S)$ functions as an infinitesimal
invariant for the map induced by $r_{H,\eta}$ (as $\eta$ varies).
For $\nu > 1$, $r_{\MF}$'s induced map serves as a
``higher infinitesimal invariant''.

\newpage

\head\S6. A map from cycle classes to 
de Rham/Mumford-Griffiths invariants
\endhead
\medbreak

In order to prove Theorem 0.4, we need cycle classes which capture cycles
with trivial Mumford-Griffiths invariants. In this section we 
construct it by using
the arithmetic cycle classes in higher extension groups in the 
category of arithmetic mixed Hodge modules. We then apply these ideas to
to a particular setting in \S7, by proving Theorem 7.2 below, which in turn implies
Theorem 0.4.

Let the notation be as in \S5.
Let $f:X\to S$ be a proper smooth morphism of quasiprojective smooth 
varieties over $k$ and let $g:S\to \Spec(k)$ be the natural morphism.
Recall the notation in \S1(iv).
For $L\in \MHM(S)$ we have the Leray spectral sequence 
$$
E_2^{p,q}=\Ext^{p}_{\MHM(\Spec(k))}(\QQ(0),R^{q}g_{\ast} L)
\Rightarrow \Ext^{p+q}_{\MHM(S)}(\QQ_S(0),L)
$$
associated to the derived functor 
$Rg_{\ast} : D^b(\MHM(S))\to D^b(\MHM(\Spec(k)))$.
Note that $\MHM(\Spec(k))$ is identified with a subcategory of $\MHS$ via 
the functor $r_{H,\eta}$ (5.4.1) with $\eta:\Spec(\CC)\to \Spec(k)$ induced 
by the natural embedding $k\hookrightarrow \CC$, and $\QQ(0)$ denotes the usual
Hodge structure. Noting that $\Ext^{p}_{\MHM(\Spec(k))}=0$ for $p\geq 2$, 
we get the short exact sequence
$$
0 \to \Ext^{1}_{\MHM(\Spec(k))}(\QQ_S(0),R^{q-1}g_{\ast} L)\to 
\Ext^{q}_{\MHM(S)}(\QQ_S(0),L) \to 
\hom_{\MHS(\Spec(k))}(\QQ(0),R^{q-1}g_{\ast} L)\to 0.
$$
Apply this to $L=R^{2r-\nu}f_{\ast}\QQ_X(r)\;\in \MHM(S)$ and we get 
natural maps
$$ 
\pi: \Ext^{\nu}_{\MHM(S)}(\QQ_S(0),R^{2r-\nu}f_{\ast}\QQ_X(r))\to 
\Lambda^{r,\nu}(X/S)
:= \hom_{\MHS}(\QQ(0),H^{\nu}(\SC,R^{2r-\nu}f_{\ast}\QQ(r))),
$$
$$ 
\ker(\pi) \to \Xi^{r,\nu}(X/S) :=
 \Ext^1_{\MHS}(\QQ(0),H^{\nu-1}(\SC,R^{2r-\nu}f_{\ast}\QQ(r))),
$$
where $H^\bullet(\SC,R^{2r-\nu}f_{\ast}\QQ(r)))$, the cohomology of $S_{an}$
with coefficients in a local system, is endowed in a canonical way with 
a mixed Hodge structure by the theory of mixed Hodge modules (\cite{MSa1}).
The composite of $\rho^{r,\nu}_{X/S}$ and $\pi$
$$ 
\lambda^{r,\nu}_{X/S}: \FB^\nu \CH^r(X/S;\QQ) \to  \Lambda^{r,\nu}(X/S) 
\tag{6.1}
$$
is given as follows. Let
$$ 
c^r_X : \CH^r(X;\QQ) \to H^{2r}(X/\CC,\QQ(r))
$$
be the cycle class map.  We  have
$c^r_X(\FB^\nu \CH^r(X/S;\QQ))\subset F_L^\nu  H^{2r}(X/\CC,\QQ(r))$,
(again, by the same argument as \cite{Sa3}, Prop. 2.1),
where $F_L^\nu$ denotes the filtration inducing the Leray spectral sequence
$$ 
E_2^{p,q}= H^p(S/\CC,R^q f_{\ast}\QQ(r))\Rightarrow
H^{p+q}(X/\CC,\QQ(r)),
$$
which degenerates at $E_2$ by Deligne's criterion (\cite{De3}).
Thus the projection to the graded quotient induces
$$ 
\FB^\nu \CH^r(X/S;\QQ) \to H^{\nu}(S/\CC,R^{2r-\nu}\rho_{\ast}\QQ(r))
$$
which is identified with $\lambda^{r,\nu}_{X/S}$ composed with 
$\Lambda^{r,\nu}(X/S)\hookrightarrow 
H^{\nu}(S/\CC,R^{2r-\nu}\rho_{\ast}\QQ(r))$, the natural inclusion. 
In particular we get the natural map
$$ 
\epsilon^{r,\nu}_{X/S}:
\FB^\nu \CH^r(X/S;\QQ)\cap \CH_{\hom}^r(X;\QQ) \to \Xi^{r,\nu}(X/S) 
\tag{6.2}
$$
for which we have 
$$\epsilon^{r,\nu}_{X/S}
\big(\FB^{\nu+1} \CH^r(X/S;\QQ) \cap \CH_{\hom}^r(X;\QQ)\big)=0.
$$

The reader can also consult \cite{Le2} for the similar construction
in terms of absolute Hodge cohomology.
\medbreak

\proclaim{Proposition 6.3} {\rm (i)} 
Assume $S=\Spec(R)$ is an affine variety. There is an injective map 
$$
\Phi^{r,\nu}_{\XS} \;:\;
\Lambda^{r,\nu}(X/S) \hookrightarrow \nabla DR^{r,\nu}(X/S)\otimes_{k}\CC,
$$
which fits into the commutative diagram
$$
\matrix
\FB^\nu\CH^r(\XS;\QQ) &@>{\phi^{r,\nu}_{\XS,\DR}}>>& 
\nabla DR^{r,\nu}(X/S)\\ 
&\\
\downarrow\rlap{$\lambda^{r,\nu}_{\XS}$}&&
{\buildrel {}_{\ \cap}\over\downarrow}\\
 \Lambda^{r,\nu}(\XS) &@>{\Phi^{r,\nu}_{\XS}}>>& 
\nabla DR^{r,\nu}(X/S)\otimes_{k}\CC \\
\endmatrix
$$
where
$\phi^{r,\nu}_{\XS,\DR}$ is defined in (3.7.1).

\medbreak\noindent
{\rm (ii)} Suppose that $S$ and $X$ are projective, and
$U \subset S$ is an open affine subvariety over $k$. Then there is a map
$$
\Phi^{r,\nu}_{X/U} \; : \;
\Lambda^{r,\nu}(X/S) \to \nabla J^{r,\nu}(X_U/U)\otimes_{k}\CC,
\quad (X_U=X\times_S U)
$$
which fits into the commutative diagram
$$
\matrix
\FB^\nu\CH^r(\XS;\QQ) &@>{\phi^{r,\nu}_{X_U/U}}>>& \nabla J^{r,\nu}(X_U/U)\\ 
&\\
\downarrow\rlap{$\lambda^{r,\nu}_{\XS}$}&&\downarrow\\
 \Lambda^{r,\nu}(\XS) &@>{\Phi^{r,\nu}_{X/U}}>>& 
\nabla J^{r,\nu}(X_U/U)\otimes_{k}\CC \\
\endmatrix
$$
\endproclaim

\demo{Proof} Part (i). Recall the notation in \S1(iv).
The complex of holomorphic forms on $S/\CC$ with
values in $R^{2r-\nu}f_{\ast}\CC$ furnishes a
resolution of $R^{2r-\nu}f_{\ast}\CC$. Since $S$ is affine, a 
variant of the Grothendieck Algebraic de Rham Theorem implies
that $H^{\nu}(\SC,R^{2r-\nu}f_{\ast}\CC)$ can be computed from 
algebraic differential forms (\cite{De2} (\S6, pp. 98-99)). Specifically 
$$
H^{\nu}(S/\CC,R^{2r-\nu}f_{\ast}\CC) \simeq 
\HH^{\nu}(\Omega_{S/k}^{\bullet}\otimes_{{\Cal O}_{S}}
\RR^{2r-\nu}f_{\ast}\Omega^{\bullet}_{\XS})\otimes_k \CC,
$$
where the latter term is computed in the Zariski topology. Since
$S$ is affine and 
$\Omega_{S/k}^{\bullet}\otimes_{{\Cal O}_{S}}
\RR^{2r-\nu}f_{\ast}\Omega^{\bullet}_{\XS}$ is coherent, 
$\HH^{\nu}(\Omega_{S/k}^{\bullet}\otimes_{{\Cal O}_{S}}\RR^{2r-\nu}
f_{\ast}\Omega^{\bullet}_{\XS})$
is computed as
$$
\frac{\ker \big(\nabla : H^{0}(S,\Omega^{\nu}_{S/k}\otimes 
\RR^{2r-\nu}f_{\ast}\Omega^{\bullet}_{\XS}) \to 
H^{0}(S,\Omega^{\nu+1}_{S/k}\otimes 
\RR^{2r-\nu}f_{\ast}\Omega^{\bullet}_{\XS})\big)}
{\nabla\big(H^{0}(S,\Omega^{\nu-1}_{S/k}\otimes 
\RR^{2r-\nu}f_{\ast}\Omega^{\bullet}_{\XS})\big)},
$$
which coincides with $\nabla DR^{r,\nu}(X/S)$ by definition (Def. 3.2).
Thus we have 
$$
\Lambda^{r,\nu}(X/S) \subset H^{\nu}(S/\CC,R^{2r-\nu}
f_{\ast}\CC) \simeq \nabla DR^{r,\nu}(X/S)\otimes_{k}\CC,
$$
which defines the injection 
$\Lambda^{r,\nu}(X/S)\hookrightarrow\nabla DR^{r,\nu}(X/S)\otimes_{k}\CC$.
\bigskip
Part (ii). From the work of Deligne (see \cite{B-Z}), the complex
$$
\big(\Omega_{S/k}^{\bullet}\otimes \RR^{2r-\nu}f_{\ast}
\Omega^{\ast\geq r-\bullet}_{X_{\CC}/S_{\CC}},\nabla\big)
$$
computes $F^{r}H^{\nu}(S/\CC,R^{2r-\nu}f_{\ast}\CC)$. This is in
the strong topology, but using the second spectral sequence
on hypercohomology  given in \cite{G-H} (p. 446):
$$
{}^{\p\p}E_{2}^{\bullet,\bullet} = 
H_{\nabla}^{\bullet}\big({H}^{\bullet}(
\Omega_{S/k}^{\bullet}\otimes \RR^{2r-\nu}f_{\ast}
\Omega^{\ast\geq r-\bullet}_{X_\CC/S_\CC})\big) \Rightarrow
 \HH^{\bullet}
\big(\Omega_{S/k}^{\bullet}\otimes \RR^{2r-\nu}f_{\ast}
\Omega^{\ast\geq r-\bullet}_{X_\CC/S_\CC}\big),
$$
and the assumed projectivity of $S$, together with Serre's work (GAGA), 
the same result holds in the Zariski topology. Specifically, we get
$$
F^{r}H^{\nu}(S/\CC,R^{2r-\nu}f_{\ast}\CC) \simeq 
\HH^{\nu}\big(\Omega_{S/k}^{\bullet}\otimes \RR^{2r-\nu}f_{\ast}
\Omega^{\ast\geq r-\bullet}_{\XS}\big)\otimes_k \CC.
$$ 
Next, by restriction there is a map
$$
\HH^{\nu}
\big(\Omega_{S/k}^{\bullet}\otimes \RR^{2r-\nu}f_{\ast}
\Omega^{\ast\geq r-\bullet}_{\XS}\big)
\to \HH^{\nu}
\big(\Omega_{U/k}^{\bullet}\otimes \RR^{2r-\nu}f_{\ast}
\Omega^{\ast\geq r-\bullet}_{X_U/U}\big)
$$
$$
\simeq  \frac{\ker \big(\nabla : 
H^{0}(U,\Omega^{\nu}_{U/k}\otimes 
\RR^{2r-\nu}f_{\ast}\Omega^{\bullet\geq r-\nu}_{X_U/U}) \to 
H^{0}(S,\Omega^{\nu+1}_{U/k}\otimes 
\RR^{2r-\nu}f_{\ast}\Omega^{\bullet\geq r-\nu-1}_{X_U/U})\big)}
{\nabla\big(H^{0}(U,\Omega^{\nu-1}_{U/k}\otimes 
\RR^{2r-\nu}f_{\ast}\Omega^{\bullet\geq r-\nu+1}_{X_U/U})\big)},
$$
where the isomorphism uses the fact that $U$ is affine.
The last term is identified with $\nabla J^{r,\nu}(X_U/U)$ by Def. 3.2.
This gives rise to the desired map.\qed
\enddemo
\medbreak

It would be nice if one could have 
Proposition 6.3 (ii) without assuming $X$ and $S$ are projective.
Indeed we will provide strong evidence in support of the following:

\proclaim{Conjecture 6.4}
Put
$$\Lambda_{\alg}^{r,\nu}(X/S)=
\Image(\FB^\nu\CH^r(\XS;\QQ) @>{\lambda^{r,\nu}_{\XS}}>> \Lambda^{r,\nu}(\XS).$$
Then, for any dense affine open $U\subset S$, there is a map
$$
\Phi^{r,\nu}_{X/U} \; :\; 
\Lambda_{\alg}^{r,\nu}(X/S) \to \nabla J^{r,\nu}(X_U/U)\otimes_{k}\CC,
\quad (X_U=X\times_S U)
$$
which fits into the commutative diagram
$$
\matrix
\FB^\nu\CH^r(\XS;\QQ) &@>{\phi^{r,\nu}_{X_U/U}}>>& \nabla J^{r,\nu}(X_U/U)\\ 
&\\
\downarrow\rlap{$\lambda^{r,\nu}_{\XS}$}&&\downarrow\\
 \Lambda_{\alg}^{r,\nu}(\XS) &@>{\Phi^{r,\nu}_{X/U}}>>& 
\nabla J^{r,\nu}(X_U/U)\otimes_{k}\CC \\
\endmatrix
%\tag{}
$$
\endproclaim

Our approach to the above conjecture involves $L^{2}$-cohomology techniques, 
where the case $\dim S = 1$ has been worked out by Zucker (see \cite{B-Z}). 
The story for $\dim S > 1$ is rather complicated, and is clarified
to some degree in the Appendix to this paper.
Indeed, the above conjecture is a consequence of other
``very reasonable'' conjectures.
\medbreak

Choose a morphism $\eta:\Spec(\CC) \to S$ whose image is the generic point
and let $X_{\eta}$ be the base change via $\eta$.
By definition $X_{\eta}$ is proper smooth over $\CC$. 
Recall $\nabla J^{r,\nu}(X_{\eta}) :=
\nabla J^{r,\nu}(X_{\eta}/\CC)$ and $\nabla DR^{r,\nu}(X_{\eta})
:= \nabla DR^{r,\nu}(X_{\eta}/\CC)$ in Definition 3.2.
%and let $\tau : \nabla J^{r,\nu}(X_{\eta})\to \nabla DR^{r,\nu}(X_{\eta})$, 
\medbreak

\proclaim{Proposition 6.5} {\rm (i)} There is a map
$$ 
\Phi^{r,\nu}_{X_{\eta}} : \Lambda^{r,\nu}(X/S) \to 
\nabla DR^{r,\nu}(X_{\eta})\otimes_{k}\CC,
$$
which fits into the commutative diagram
$$
\matrix
\FB^\nu \CH^r(X/S;\QQ) &@>{\lambda^{r,\nu}_{X/S}}>>& 
\Lambda^{r,\nu}(X/S) \\
&\\
\downarrow\rlap{$\eta^{\ast}$}&&\downarrow\rlap{$\Phi^{r,\nu}_{X_{\eta}}$}\\
&\\
\FB^\nu\CH^r(X_{\eta};\QQ) &@>{\phi^{r,\nu}_{X_{\eta}}}>>
\nabla J^{r,\nu}(X_{\eta})\otimes_{k}\CC @>>>& 
\nabla DR^{r,\nu}(X_{\eta})\otimes_{k}\CC \\
\endmatrix
%\tag{6.2}
$$
In particular we have
$$\phi^{r,\nu}_{X_{\eta},DR}\circ\eta^{\ast}(\FB^\nu 
\CH^r(X/S;\QQ)\cap \CH_{\hom}^r(X;\QQ))=0.$$

\medbreak\noindent
{\rm (ii)} Suppose that $X$ and $S$ are projective over $k$ or that
Conjecture 6.4 is true. Then there is a morphism
$$
\Lambda_{\alg}^{r,\nu}(X/S) \to \nabla J^{r,\nu}(X_{\eta})\otimes_{k}\CC
$$ 
which is compatible with the diagram in 6.5. In particular, we have
$$\phi^{r,\nu}_{X_{\eta}}\circ\eta^{\ast}(
\FB^\nu \CH^r(X/S;\QQ)\cap \CH_{\hom}^r(X;\QQ))=0.$$
\endproclaim

\medbreak\noindent
\underbar{Remarks 6.7}. 
The construction of the map
$\Lambda^{r,\nu}(X/S) \to \nabla J^{r,\nu}(X_{\eta})\otimes_k\CC$, also
appears in \cite{Ke1} in the special case
where $X = Y\times S$ with $Y,\ S\in \C_{k}$.
In this case, for an affine open subscheme $U\subset S$, we have
$$
\align
\Lambda^{r,\nu}(X\times_S U/U) = 
\hom_{\MHS}(\QQ(0),H^{\nu}(U/\CC,\QQ)\otimes H^{2r-\nu}(Y/\CC,\QQ))
&\hookrightarrow \HH^{\nu}(\Omega_{U/k}^{\bullet})\otimes 
F^{r-\nu}H^{2r-\nu}_{\DR}(Y/k)\otimes_k\CC \\
= H^{\nu}\big(&\Gamma(U,\Omega_{U/k}^{\bullet})\big)
\otimes F^{r-\nu}H^{2r-\nu}_{\DR}(Y/k)\otimes_k\CC,
\endalign
$$
via the comparison isomorphism
$H^{\bullet}_{\DR}(Y/k)\otimes_k\CC\simeq H^{\bullet}(Y/\CC,\CC)$
and Grothendieck's algebraic de Rham theorem.  One then applies the
morphism
$$
H^{\nu}\big(\Gamma(U,\Omega_{U/k}^{\bullet})\big)
\otimes F^{r-\nu}H^{2r-\nu}_{\DR}(Y/k)\otimes_k\CC \to
H^{\nu}_{\nabla}\big(\Omega_{U/k}^{\bullet}
\otimes F^{r-\bullet}H^{2r-\nu}_{\DR}(Y/k)\big)\otimes_k\CC.
$$

\head\S7. Detecting nonzero classes with trivial de Rham/Mumford 
invariant\endhead

In this section we take the base field $k=\ol{\QQ}$.
Take $X \in {\Cal C}$ with $\dim X = d$. 
For $\Gamma \in \C_{\ol{\QQ}}$ with $\dim \Gamma = 1$ and $S\in \C$ and 
$\xi\in \CH^{r}(\Gamma/\CC\times S\times X;\QQ)$, 
consider the cycle induced map
$$
[\xi]_{\ast} : H^{1}(\Gamma/\CC,\QQ)\otimes H^{2d_S-\nu+1}(S,\QQ)
\to H^{2r-\nu}(X,\QQ).\tag{7.0}
$$
where $d_S=\dim S$. It is easy to see that 
$$
[\xi]_{\ast}\big( H^{1}(\Gamma/\CC,\QQ)\otimes H^{2d_S-\nu+1}(S,\QQ)
\big) \subset H^{2r-\nu}(X,\QQ)\cap F^{r-\nu}H^{2r-\nu}(X,\CC).
$$
Now introduce:

\proclaim{Definition 7.1} Write $H^{r-\nu,r}(X)=H^r(X,\Omega^{r-\nu}_{X/\CC})$.
\medbreak\noindent
{\rm (i)} Define 
$H_{\ol{\QQ}}^{r-\nu,r}(X)\subset H^{r-\nu,r}(X)$
as the $\CC$-vector space generated by the Hodge projected images of 
$[\xi]_{\ast}$ in 7.0, over all $\Gamma\in \C_{\ol{\QQ}}$ of dimension $1$, 
and all $S\in \C$, and all $\xi\in \CH^{r}(\Gamma/\CC\times S\times X;\QQ)$.
\medbreak\noindent
{\rm (ii)} 
For $X_o\in \C_{\ol{\QQ}}$, we introduce a variant 
$H^{r-\nu,r}_{\ol{\QQ}}(X_o/\ol{\QQ})\subset H^{r-\nu,r}(X)$
($X=X_o\times_{\ol{\QQ}}\CC$)
of the above definition, where we only consider $S\in \C_{\ol{\QQ}}$
and where we allow $\xi\in \CH^{r}((\Gamma\times S\times X)/\CC;\QQ)$.
\endproclaim
\medbreak

For $X\in \C$ we recall the maps (see (3.7.2) and (3.7.3))
$$
\phi^{r,\nu}_{X,DR}: \GrFB^{\nu}\CH^{r}(X;\QQ)\to \nabla DR^{r,\nu}(X)
\quad\text{ and }\quad
\phi^{r,\nu}_{X}: \GrFB^{\nu}\CH^{r}(X;\QQ)\to \nabla J^{r,\nu}(X).
$$
The main result is the following:

\proclaim{Theorem 7.2} Assume $C(X)$ (cf. \S1 (v)).
\medbreak\noindent
{\rm (i)} If $\nu \geq 2$ and if $H_{\ol{\QQ}}^{r-\nu,r}(X)\ne 0$, then
there are an uncountable number of classes in $\ker(\phi^{r,\nu}_{X,DR})$.
\medbreak\noindent
{\rm (ii)} Let 
$X=X_o\times_{\ol{\QQ}}\CC$ with $X_o\in \C_{\ol{\QQ}}$.
If $\nu \geq 2$ and if $H_{\ol{\QQ}}^{r-\nu,r}(X_o/\ol{\QQ})\ne 0$, then
there are an uncountable number of classes in 
$\ker(\phi^{r,\nu}_{X})$.
\medbreak\noindent
{\rm (iii)} Assume Conjecture 6.4 is true.
If $\nu \geq 2$ and if $H_{\ol{\QQ}}^{r-\nu,r}(X)\ne 0$, then
there are an uncountable number of classes in $\ker(\phi^{r,\nu}_{X})$.
\endproclaim

\demo{Proof of Theorem} Since the proofs of (ii) and (iii) are similar to
and simpler than (i), we will just prove (i).
By assumption, one can find $\Gamma/{\ol{\QQ}}$, $S$ and 
a cycle $\xi \in \CH^{r}(\Gamma/\CC\times S\times X)$
such that the  map
$$
\xi_{\ast} : H^{1}(\Gamma/\CC,\CC) \to 
H^{\nu-1}(S,\CC)\otimes \big(H^{2r-\nu}(X,\CC)/
F^{r-\nu+1}H^{2r-\nu}(X,\CC)\big),
\tag{7.2.1}
$$
is nonzero. Note that $\nu\geq 2$ implies $\dim S \geq 1$,
which will be exploited later.

In what follows, for a $\ZZ$- or $\QQ$-module $H$ of finite rank endowed 
with a mixed Hodge structure, we denote
$$ 
J(H)=W_0 H\otimes\CC /\big(F^0 W_0H\otimes\CC + W_0 H\big).
$$
It is naturally endowed with a structure of a complex Lie group
and for a morphism $H\to H'$, the induced map
$J(H) \to J(H')$ is a homomorphism of complex Lie groups.
We note that
$$ 
J^r(S\times X) = J(H^{2r-1}(S\times X,\ZZ)(r))
$$
and the canonical isomorphism due to Carlson for 
a $\QQ$-mixed Hodge structure $H$
$$
\Ext^1_{\MHS}(\QQ(0),H) \simeq J(H).
$$
Consider the composite map
$$
\Phi_{\xi} : \CH^{1}_{\deg 0}(\Gamma/{\ol{\QQ}}) \to 
J^{1}(\Gamma/\CC) @>{\xi_{\ast}}>>
J^{r}(S\times X) \to 
J(H^{\nu-1}(S,\ZZ)\otimes H^{2r-\nu}(X,\ZZ)(r)),
$$
where the last map is induced by the projection to the K\"unneth component.
For dense (Zariski) open $U\subset S$ let
$$ 
\Phi_{\xi,U} : \CH^{1}_{\deg 0}(\Gamma/{\ol{\QQ}}) \to 
J(H^{\nu-1}(U,\ZZ)\otimes H^{2r-\nu}(X,\ZZ)(r))
$$
be the composite of $\Phi_{\xi}$ with the natural map
$$ 
J(H^{\nu-1}(S,\ZZ)\otimes H^{2r-\nu}(X,\ZZ)(r))\to
J(H^{\nu-1}(U,\ZZ)\otimes H^{2r-\nu}(X,\ZZ)(r)).
$$
A key to the proof of Theorem 7.2 is the following:

\proclaim{Lemma 7.3}
Under the above assumption, 
$$
\lim_{\buildrel \to\over {U\subset S}}\text{\rm Image}(\Phi_{\xi,U}),
$$ 
has infinite rank.
\endproclaim
\medbreak

We now finish the proof of Theorem 7.2 assuming Lemma 7.3. 
We choose spreads of $X$ and $S$: Choose $M$ smooth over $\ol{\QQ}$ and 
proper smooth morphisms $\X\to M$ and ${\Cal S}\to M$ such that
for a suitable $\eta:\Spec(\CC)\to M$ whose image is the generic point, we have
$\X\times_M\Spec(\CC)\simeq X$ and 
${\Cal S}\times_M\Spec(\CC)\simeq S$ by the base change via $\eta$.
Correspondingly we have 
$\tilde{\xi}\in \CH^r(\Gamma\times \X\times_M {\Cal S};\QQ)$
whose image under
$$ 
\eta^{\ast}:\CH^r(\Gamma\times \X\times_M {\Cal S};\QQ)\to 
\CH^r(\Gamma\times X\times S;\QQ)
$$
coincides with $\xi$.
By \cite{Ja}, 7.2 (see also \cite{Sa2}, 5-1) one can lift the 
the K\"unneth components of the cohomology classes of the diagonal 
class $[\Delta_X]\in H^{2d}(X\times X,\QQ)$ 
(which are algebraic by $C(X)$) to 
$$
\tilde{\Delta}_{\X}(i,j)\in \CH^d(\X\times_M \X;\QQ)/
\underset{\mu\geq 1}\to{\bigcap} \FB^\mu \CH^d(\X\times_M \X/M;\QQ)
\quad (i+j=2d)
$$
such that their restrictions to the fiber $X_s$ for every $s:\Spec(\CC)\to M$
$$
\tilde{\Delta}_{X_s}(i,j)\in \CH^d(X_s\times X_s;\QQ)/
\underset{\mu\geq 1}\to{\bigcap} \FB^\mu \CH^d(X_s\times X_s;\QQ) 
\quad (i+j=2d)
$$
give a Chow-K\"unneth decomposition of the diagonal.
Then we can let $\tilde{\Delta}_{\X}(2d-2r+\nu,2r-\nu)$ act on  $\tilde{\xi}$ 
(and replace $\tilde{\xi}$ by the result) without changing (7.2.1). 
This will ensure that $\tilde{\xi}$ induces
$$ 
\Phi_{\tilde{\xi}} : \CH^{1}_{\deg 0}(\Gamma/{\ol{\QQ}}) \to 
F_{B}^{\nu}\CH^{r}(\X\times_M {\Cal S}/{\Cal S};\QQ)
\to \GrFB^{\nu}\CH^{r}(\X\times_M {\Cal S}/{\Cal S};\QQ)
$$
and its image is contained in $\CH_{\hom}^r(\X\times_M {\Cal S};\QQ)$. 
By the construction in \S6, this means
$$
\Image(\Phi_{\tilde{\xi}}) \subset \ker\big(
\GrFB^{\nu}\CH^{r}(\X\times_M {\Cal S}/{\Cal S};\QQ)
@>{\phi^{r,\nu}_{\X\times_M {\Cal S}/{\Cal S},DR}}>>
\nabla DR^{r,\nu}(\X\times_M {\Cal S}/{\Cal S})\big).
\tag{*}
$$
By construction, $\Phi_{\xi}$ coincides with the composite of 
$\Phi_{\tilde{\xi}}$ and the map (cf. (6.2))
$$ 
\align
\FB^{\nu}\CH^{r}(\X\times_M {\Cal S}/{\Cal S};\QQ)\cap\CH_{\hom}^r(\X\times_M 
{\Cal S};\QQ)
@>{\epsilon^{r,\nu}_{\X\times_M {\Cal S}/{\Cal S}}}>>
&\Xi^{r,\nu}(\X\times_M {\Cal S}/{\Cal S})\\
:=&\Ext^1_{\MHS}(\QQ(0),H^{\nu-1}({\Cal S}/\CC,R^{2r-\nu}f_{\ast}\QQ(r))),\\
\endalign
$$
where $f:\X\times_M{\Cal S} \to {\Cal S}$ is the projection, and the map
$$ 
\Xi^{r,\nu}(\X\times_M {\Cal S}/{\Cal S}) \to 
\Ext^1_{\MHS}(\QQ(0),H^{\nu-1}(S,\QQ)\otimes H^{2r-\nu}(X,\QQ)(r))=
 J(H^{\nu-1}(S,\QQ)\otimes H^{2r-\nu}(X,\QQ)(r))
$$
which is the restriction via $S \hookrightarrow {\Cal S}/\CC$ 
induced by $\eta$ 
(Recall that $\eta:\Spec(\CC)\to M$ gives rise to 
$\eta_{\CC}:\Spec(\CC)\to M/\CC$ and 
$S\simeq {\Cal S}/\CC\times_{M/\CC} \Spec(\CC)$ by the base change via 
$\eta_{\CC}$).
Noting that 
$\epsilon^{r,\nu}_{\X\times_M {\Cal S}/{\Cal S}}$ annihilates
$\FB^{\nu+1}\CH^{r}(\X\times_M {\Cal S}/{\Cal S};\QQ))$,
Lemma 7.3 implies that the image of 
$$
\CH^{1}_{\deg 0}(\Gamma/{\ol{\QQ}}) @>{\Phi_{\tilde{\xi}}}>>
\GrFB^{\nu}\CH^{r}(\X\times_M {\Cal S}/{\Cal S};\QQ) @>{j^{\ast}}>> 
\GrFB^{\nu}\CH^{r}(\X\times_M {\Cal U}/{\Cal U};\QQ) 
$$
has infinite rank in the limit over any dense open $j:{\Cal U}\hookrightarrow {\Cal S}$.
Identify a point $p\in S$ with a morphism $p:\Spec(\CC)\to S$ and let
$\tilde{p} : \Spec(\CC)\to {\Cal S}$ be the composite of $p$ and 
$S \to {\Cal S}$ induced by $\eta$. Then the composite of 
$\tilde{p}$ and ${\Cal S}\to M$ is $\eta$ so that we have
$X\simeq \X\times_M {\Cal S}\times_{\Cal S} \Spec(\CC)$ by the base change via 
$\tilde{p}$. If the image of $\tilde{p}$ is the generic point, it induces
the map 
$$
\tilde{p}^{\ast}: \lim_{\buildrel\longrightarrow\over{{\Cal U}\subset{\Cal S}}}
\GrFB^\nu\CH^{r}(\X\times_M {\Cal U};\QQ)\to \GrFB^\nu\CH^r(X;\QQ),
$$
where ${\Cal U}$ ranges over the dense open subschemes of ${\Cal S}$, and it is 
injective by [Sa3], Rem.(2-1)(3). Therefore the composite map
$$ 
\Phi_{\tilde{\xi},p}:
\CH^{1}_{\deg 0}(\Gamma/{\ol{\QQ}}) @>{\Phi_{\tilde{\xi}}}>> 
\GrFB^{\nu}\CH^{r}(\X\times_M {\Cal S}/{\Cal S};\QQ) 
@>{\tilde{p}^{\ast}}>> \GrFB^{\nu}\CH^{r}(X;\QQ)
$$
has the image of infinite rank and it is contained in 
$\ker(\phi^{r,\nu}_{X,DR})$ by $(*)$ and Proposition 6.5.
Noting that $\Phi_{\tilde{\xi},p}$ coincides with 
$$ 
\CH^{1}_{\deg 0}(\Gamma/{\ol{\QQ}}) \to 
\CH^{1}_{\deg 0}(\Gamma/{\ol{\QQ}})\otimes \CH_0(S;\QQ) 
@>{\xi^{\ast}}>> \GrFB^\nu\CH^{r}(X;\QQ)
$$
where the first map sends $\alpha\in \CH^{1}_{\deg 0}(\Gamma/{\ol{\QQ}})$
to $\alpha\otimes [p]$, we have thus shown that the image of
$$
\xi_{\ast} : \CH^{1}_{\deg 0}(\Gamma/{\ol{\QQ}})\otimes \CH_{0}(S;\QQ) \to 
\GrFB^{\nu}\CH^{r}(X;\QQ),
$$
is nontrivial and lies in $\ker(\phi^{r,\nu}_{X,DR})$.
Thus for suitable $P-Q \in \CH_{\deg 0}^{1}(\Gamma/{\ol{\QQ}})$,
we have a nontrivial cycle induced map
$$
\CH_{0}(S;\QQ)\to\ker(\phi^{r,\nu}_{X,DR})\subset \GrFB^{\nu}\CH^{r}(X;\QQ).
$$
This induces a nontrivial composite map 
$$
\lambda : S \to \CH_{0}(S/{\CC};\QQ) \to \ker(\phi^{r,\nu}_{X,DR})
\hookrightarrow \GrFB^{\nu}\CH^{r}(X;\QQ), 
%\hookrightarrow 
%\CH^{r}(X;\QQ)/\underset_{\nu>0}\to{\bigcap} \FB^\nu\CH^{r}(X;\QQ),
%\tag{7.8}
$$
whose image by an easy argument using $\dim S\geq 1$, is more than a point.
%Just introduce $[\Delta_{S}-\Delta_{S}(0,2d_{S})]$, which is algebraic.
The standard arguments in the theory of Chow varieties
imply that the fibers of $\lambda$ are $c$-closed
(countably closed). I.e. the fibers of $\lambda$ are countable
unions of {\it proper} subvarieties of $S$. This
is the import of the theory in \cite{R1-2}, (also 
cf. \cite{Sch1}), and it basically hinges on Lemma 7.4 below.
Thus if the image of $\lambda$ were countable, then $S$ would be a countable 
union of proper subvarieties. This is impossible by Baire's Theorem, 
which completes the proof of Theorem 7.2.\qed

\proclaim{Lemma 7.4}
Let $X\in \C$ and $Y$ be a projective variety over $\CC$, and assume
given a cycle induced map
$$
\kappa : Y \to \CH^{r}(X;\QQ).
$$
Then $\kappa^{-1}\big(\FB^{\nu}\CH^{r}(X;\QQ)\big)$ is a $c$-closed
subset of $Y$. Furthermore, if $\kappa (Y) 
\subset \FB^{\nu}\CH^{r}(X;\QQ)$, then the fibers of the induced map
$Y \to \GrFB^{\nu}\CH^{r}(X;\QQ)$
are $c$-closed.
\endproclaim
\demo{Proof} Omitted.\qed
\enddemo

\enddemo
\medbreak

Now we prove Lemma 7.3. 

\proclaim{Sublemma 7.3.1} Consider the map
$$
\xi_{\ast} : H^{1}(\Gamma/\CC,\QQ)
\to H^{\nu-1}(S,\QQ)\otimes H^{2r-\nu}(X,\QQ) 
$$
induced by $\xi$ (as in eqn. (7.2.1)). Then its image is not contained in 
$$
N^1_H H^{\nu-1}(S,\QQ)\otimes H^{2r-\nu}(X,\QQ) +
H^{\nu-1}(S,\QQ)\otimes N^{r-\nu+1}_H H^{2r-\nu}(X,\QQ),
$$
where we recall $N_{H}^{p}H^{\bullet}(X,\QQ)$ is the largest
subHodge structure of $H^{\bullet}(X,\QQ)$
lying in $F^{p}H^{\bullet}(X,\CC) \cap 
H^{\bullet}(X,\QQ)$. 
\endproclaim
\demo{Proof of Claim} Obvious.\qed
\enddemo

\proclaim{Sublemma 7.3.2} 
The image of 
$$J^1(\Gamma/{\CC}) @>{\xi_{\ast}}>> J^{r}(S\times X) \to 
J(H^{\nu-1}(S,\ZZ)\otimes H^{2r-\nu}(X,\ZZ)(r))) \to
J(H^{\nu-1}(U,\ZZ)\otimes H^{2r-\nu}(X,\ZZ)(r)))
$$
is nonzero. 
\endproclaim

\demo{Proof} Recall that 
$$
J(H^{\nu-1}(S,\ZZ)\otimes H^{2r-\nu}(X,\ZZ)(r))=V/L
$$
with
$$
V:=\big(H^{\nu-1}(S,\CC)\otimes H^{2r-\nu}(X,\CC)\big)/
F^r\big(H^{\nu-1}(S,\CC)\otimes H^{2r-\nu}(X,\CC)\big)
$$
and $L\subset V$ a lattice. We observe that
$$ 
F^r\big(H^{\nu-1}(S,\CC)\otimes H^{2r-\nu}(X,\CC)\big)\subset
H^{\nu-1}(S,\CC)\otimes F^{r-\nu+1}H^{2r-\nu}(X,\CC).
$$
By the assumption (cf. (7.2.1)) the image of 
$J^1(\Gamma/{\CC}) @>{\xi_{\ast}}>> J^{r}(S\times X) \to 
J(H^{\nu-1}(S,\ZZ)\otimes H^{2r-\nu}(X,\ZZ)(r)))
$
is nonzero. By Sublemma 7.3.1 it now suffices to show that the kernel of 
$$ 
J(H^{\nu-1}(S,\ZZ)\otimes H^{2r-\nu}(X,\ZZ)(r))\to
J(H^{\nu-1}(U,\ZZ)\otimes H^{2r-\nu}(X,\ZZ)(r)).
$$
is contained in the image of
$$
J(N_H^{1}H^{\nu-1}(S,\QQ)\otimes H^{2r-\nu}(X,\QQ)) + 
J(H^{\nu-1}(S,\QQ)\otimes N_{H}^{r+1-\nu}H^{2r-\nu}(X,\QQ)).
$$
Indeed one has a surjective map
$$
H^{\nu-1}(S,\QQ)\otimes H^{2r-\nu}(X,\QQ)(r) \twoheadrightarrow
W_{\nu-1}H^{\nu-1}(U,\QQ)\otimes H^{2r-\nu}(X,\QQ)(r),
$$
whose kernel is contained in 
$N_H^{1}H^{\nu-1}(S,\QQ)\otimes H^{2r-\nu}(X,\QQ)(r)$. 
Correspondingly the kernel of
$$ 
J(H^{\nu-1}(S,\QQ)\otimes H^{2r-\nu}(X,\QQ)(r)) \to 
J(W_{\nu-1}H^{\nu-1}(U,\QQ)\otimes H^{2r-\nu}(X,\QQ)(r))
$$
is contained in the image of 
$J(N_H^{1}H^{\nu-1}(S,\QQ)\otimes H^{2r-\nu}(X,\QQ)).$
Next, there is a short exact sequence
$$
0 \to W_{\nu-1}H^{\nu-1}(U,\QQ)\otimes H^{2r-\nu}(X,\QQ)(r) \to
W_{0}\big(H^{\nu-1}(U,\QQ)\otimes H^{2r-\nu}(X,\QQ)(r)\big)
$$
$$
\to Gr_{W}^{0}\big(H^{\nu-1}(U,\QQ)\otimes H^{2r-\nu}(X,\QQ)(r)\big) \to 0.
$$
Taking global $\Ext_{\MHS}$, it gives us the exact sequence
$$
\hom_{\MHS}\big(\QQ(0),Gr_{W}^{0}(
H^{\nu-1}(U,\QQ)\otimes H^{2r-\nu}(X,\QQ)(r))
\big) \to 
$$
$$
J\big(W_{\nu-1}H^{\nu-1}(U,\QQ)\otimes H^{2r-\nu}(X,\QQ)(r)\big)
\to J\big(H^{\nu-1}(U,\QQ)\otimes H^{2r-\nu}(X,\QQ)(r)\big).
$$
Hence it suffices to show that the image of the first map is contained in
$$
J\big(W_{\nu-1}H^{\nu-1}(U,\QQ)\otimes 
N_{H}^{r+1-\nu}H^{2r-\nu}(X,\QQ)(r)\big).
$$
We have a short exact sequence
$$
0 \to W_{\nu-1}H^{\nu-1}(U,\QQ)\otimes N_{H}^{r+1-\nu}H^{2r-1}(X,\QQ)(r) \to
W_{0}\big(H^{\nu-1}(U,\QQ)\otimes N_{H}^{r+1-\nu}H^{2r-1}(X,\QQ)(r)\big)
$$
$$
\to Gr_{W}^{0}\big(H^{\nu-1}(U,\QQ)\otimes N_{H}^{r+1-\nu}
H^{2r-\nu}(X,\QQ)(r)\big) \to 0.
$$
Moreover, by using the fact $F^p H^{\nu-1}(U,\CC)/F^{p+1}=0$ for $p>\nu-1$, 
one sees that the natural inclusion:
$$
{}_{\hom_{\MHS}\big(\QQ(0),Gr_{W}^{0}(H^{\nu-1}(U,\QQ)\otimes N_{H}^{r+1-\nu}
H^{2r-\nu}(X,\QQ)(r))\big)
\quad\subset\quad
\hom_{\MHS}\big(\QQ(0),Gr_{W}^{0}(H^{\nu-1}(U,\QQ)\otimes H^{2r-\nu}(X,\QQ)(r))
\big),}
$$
is an equality. 
Finally the desired assertion follows from the commutative diagram:
$$
\matrix
{}_{\hom_{\MHS}\big(\QQ(0),Gr_{W}^{0}(H^{\nu-1}(U,\QQ)\otimes 
N_{H}^{r+1-\nu}H^{2r-\nu}(X,\QQ)(r))
\big)}&{}_{\to}& {}_{J\big(W_{\nu-1}H^{\nu-1}(U,\QQ)\otimes 
N_{H}^{r+1-\nu}H^{2r-\nu}(X,\QQ)(r)\big)}\\
&\\
||&&\downarrow\\
&\\
{}^{\hom_{\MHS}\big(\QQ(0),Gr_{W}^{0}(H^{\nu-1}
(U,\QQ)\otimes H^{2r-\nu}(X,\QQ)(r))
\big)}&{}^{\to}&{}^{J\big(
W_{\nu-1}H^{\nu-1}(U,\QQ)\otimes H^{2r-\nu}(X,\QQ)(r)\big)}\\
\endmatrix
$$\qed
\enddemo

In order to complete the proof of Lemma 7.3 we need the following key result.

\proclaim{Lemma 7.5} Let $W$ be an Abelian variety over ${\ol{\QQ}}$, and
suppose that $W/\CC$ contains an Abelian subvariety $V$ over $\CC$. 
Then $W$ contains an Abelian subvariety over $\ol{\QQ}$ of dimension $=\dim V$.
\endproclaim

\demo{Proof} Construct a $\ol{\QQ}$-spread
$$
\matrix {\Cal V}&\hookrightarrow&M\times W\\
&\\
\quad\searrow&&\swarrow\quad\\
&\\
&M
\endmatrix
$$
where $M,\ {\Cal V}$ are defined over $\ol{\QQ}$
and there is $\eta:\Spec(\CC) \to M$ whose image is the generic point of $M$,
such that ${\Cal V} \times_M \Spec(\CC) \simeq V$. Let
$s\in M(\ol{\QQ})$ be a general choice of closed point.
Then the desired Abelian subvariety of $W$ is given
by ${\Cal V}_{s}$.\qed
\enddemo

The following is probably well-known, but is an immediate consequence of 
Lemma 7.5.

\proclaim{Corollary 7.6} $W/\ol{\QQ}$ is simple
$\Leftrightarrow W\times_{\ol{\QQ}} \CC$ is simple.
\endproclaim

Now we complete the proof of Lemma 7.3. 
Notice that $\CH^{1}_{\deg 0}(\Gamma/{\ol{\QQ}})$ 
has the structure of an Abelian variety over $\ol{\QQ}$. 
By Sublemma 7.3.2 and Corollary 7.6 together with Poincar\'e's complete 
reducibility theorem, there exists a simple Abelian subvariety 
$A/{\ol{\QQ}}\subset \CH^{1}_{\deg 0}(\Gamma/{\ol{\QQ}})$, 
for which the composite
$$
A/{\CC} \hookrightarrow \CH^{1}_{\deg 0}(\Gamma/{\CC})\simeq
J^1(\Gamma/{\CC})  \to 
J(H^{\nu-1}(U,\ZZ)\otimes H^{2r-\nu}(X,\ZZ)(r))
$$
is an isogeny onto its image. Now Lemma 7.3 follows from the fact that 
$A(\ol{\QQ})$ has infinite rank (\cite{F-G}).
\bigskip

\noindent
\underbar{Example 7.7}. Let $X = C\times C$, where $C$ is a smooth 
projective curve defined over $\ol{\QQ}$. In this case $C(X)$ trivially holds. 
Let 
$$[\Delta_{X}] = \sum_{p+q=4}[\Delta_{X}(p,q)]\in
\bigoplus_{p+q=4} H^2(X\times X,\QQ)
\quad 
\text{ with } \Delta_{X}(p,q)\in \CH^2(X\times X;\QQ),
$$
be the K\"unneth decomposition of the diagonal $X\hookrightarrow X\times X$. 
Using the notation in Theorem 7.2, 
we take $\Gamma = S = C$ and $r=\nu = 2$ and $\xi=\Delta_{X}(2,2)$.
Let $\eta$ be the generic point of $C$ and take $P,\ Q\in C(\ol{\QQ})$.
Noting 
$$
\Delta_{X,\ast} = \Delta_{X}(0,4)_{\ast} + \Delta_{X}(1,3)_{\ast} + 
\Delta_{X}(2,2)_{\ast}\quad \text{ on } \CH^{2}(X;\QQ),
$$
we get
$$
\xi_{\ast}\big((P-Q)\times \eta\big) = (P-Q)\times \eta
- \Delta_{X}(1,3)_{\ast}\big((P-Q)\times \eta\big).
$$
Next, we expand $\Delta_{X}(1,3)$ in terms of the diagonal
$\Delta_{C}$ of $C\hookrightarrow C\times C$. Fix a rational point $O\in
C(\ol{\QQ})$. Let $T : C\times C\times C\times C \to 
C\times C\times C\times C$ be the map given by 
$T(t_{1},t_{2},t_{3},t_{4}) = (t_{1},t_{3},t_{2},t_{4})$.
Then 
$$
\Delta_{X}(1,3)_{\ast}\big((P-Q)\times \eta\big) = 
\big[T^{\ast}\big(C\times O\times \Delta_{C} + \Delta_{C}\times 
C\times O\big)\big]_{\ast}\big((P-Q)\times \eta\big)
$$
$$
= \big[T^{\ast}\big(\Delta_{C}\times 
C\times O\big)\big]_{\ast}\big((P-Q)\times \eta\big)
= (P-Q)\times O.
$$
Thus we get
$$
\xi_{\ast}\big((P-Q)\times \eta\big)= (P-Q)\times (\eta - O),
$$
which is the cycle that was studied earlier by A. Rosenschoen and M. 
Saito \cite{RS}, and later by Matt Kerr \cite{Ke1}. It should
be pointed out that their results about a given $0$-cycle
satisfying a certain condition, 
are more specific than what is presented here. 
\bigskip
\newpage

\head\S8. A class of examples with trivial Mumford-Griffiths invariant\endhead

We first recall that if $V,\ W\in \C_k$, then (\cite{K1})
$$
B(V) \ + \ B(W) \Rightarrow B(V\times W).
$$

Another fact that we use is the following. Let $m=\dim W$,
and assume that $j: V \hookrightarrow W$ is of dimension $\ell$,
being the cutout of $W$ by $m-\ell$ hyperplane sections. By
the Lefschetz theorem, $j^{\ast} : H^{\ell}(W)\hookrightarrow
H^{\ell}(V)$ is injective.  Then
$B(W) \Rightarrow B(V)$ and the surjective left inverse
$$
\big(j^{\ast}\big)^{-1}: H^{\ell}(V)\to H^{\ell}(W),
$$
is cycle induced. Indeed, one checks that
$$
j_{\ast}\circ j^{\ast} = L_{W}^{m-\ell} : H^{\ell}(W)
\ {\buildrel \sim\over \to}\ H^{2m-\ell}(W),
$$
where $L_{W}$ is the cup product with the polarizing class.
If $\Lambda_{W}^{m-\ell}$ is the inverse to
$L_{W}^{m-\ell}$, then we obtain 
$\big(j^{\ast}\big)^{-1}= \Lambda_{W}^{m-\ell}\circ j_{\ast}$ from the formula
$$
\big(j^{\ast}\big)^{-1}\circ j^{\ast} = 
\Lambda_{W}^{m-\ell}\circ j_{\ast}\circ j^{\ast}
= \Lambda_{W}^{m-\ell}\circ L_{W}^{m-\ell} = \text{\rm Identity}.
$$

The following diagram makes the above construction transparent: 
$$
\matrix &&H^{\ell}(W)&--\\
&&&\hskip.2in\big|\\
&{}^{j^{\ast}}\swarrow&\ \wr\downarrow L_{W}^{m-\ell}&\hskip.2in\big|\\
&&&\hskip.2in\big|\\
H^{\ell}(V)&{\buildrel j_{\ast}\over\to}&H^{2m-\ell}(W)&
\hskip.2in\downarrow& \text{\rm Identity}\\
&&&\hskip.2in\big|\\
&{}_{(j^{\ast})^{-1}}\searrow&\ \wr\downarrow \Lambda_{W}^{m-\ell}\\
&&&\hskip.2in\big|\\
&&H^{\ell}(W)&\longleftarrow
\endmatrix
$$

\proclaim{Lemma 8.1} {\rm (i)}
Let $X,\; Y \in \C$ and $\Gamma\in \C_{\ol{\QQ}}$
with $\dim X = d$, $\dim \Gamma = 1$, $\dim Y = d-1$. 
Assume that  $B(Y)$ and $C(X)$ hold. Further,
suppose that we are given a dominating morphism
$$
\lambda\;:\; \Gamma_\CC \times Y \to X,
$$
such that the induced map
$$
H^{0,1}(\Gamma/\CC)\otimes H^{0,\nu-1}(Y) \to H^{0,\nu}(X),
$$
is nonzero for some integer $\nu \geq 2$. Then
$H_{\ol{\QQ}}^{r-\nu,r}(X) \ne 0$ for $\nu\leq r\leq d$.
\medbreak\noindent
{\rm (ii)}
Further, if $X=X_o\times_{\ol{\QQ}}\CC$ and $Y=Y_o\times_{\ol{\QQ}}\CC$
with $X_o,\; Y_o\in \C_{\ol{\QQ}}$, then 
$H_{\ol{\QQ}}^{r-\nu,r}(X_o/\QQ) \ne 0$ for $\nu\leq r\leq d$.
\endproclaim
\demo{Proof}
We only prove (i). The proof of (ii) is similar.
The hard Lefschetz isomorphism implies
that
$$
L_{X}^{r-\nu}: H^{0,\nu}(X) \hookrightarrow H^{r-\nu,r}(X),
$$
is injective.
Now let $S$ be a smooth cutout
of $Y$ by $d-\nu$ hyperplane sections, where $\dim S = \nu-1$.
By what we said in the beginning of this section, 
there is a nontrivial cycle induced map
$$
H^{0,1}(\Gamma/\CC)\otimes H^{0,\nu-1}(S) \to H^{r-\nu,r}(X).
$$
More explicitly, by the above discussion, the
cycle inducing this map is given as follows:
Let $\Delta_{\Gamma}$ be the diagonal of $\Gamma$
in $\Gamma\times \Gamma$. If $j : S\hookrightarrow Y$
is the inclusion, and if we identify a map $T_{\ast}$
with its associated correspondence $\{T_{\ast}\}$, 
then the cycle is
$$
\xi := L_{X}^{r-\nu}\circ \{\lambda_{\ast}\}\circ \biggl(
\biggl\langle \text{\rm Pr}_{13}^{\ast}(\Delta_{\Gamma_\CC})\bullet \text{\rm 
Pr}_{24}^{\ast}(\Lambda_{Y}^{d-\nu}\circ \{j_{\ast}\})\biggr\rangle_{\Gamma
\times S\times \Gamma\times Y}\biggr) \in \CH^{r}(\Gamma/\CC\times 
S\times X;\QQ).
$$
This gives a nonzero element of $H_{\ol{\QQ}}^{r-\nu,r}(X)$.\qed
\enddemo
\medbreak

By Theorem 7.2, we arrive at the following:

\proclaim{Corollary 8.2} 
{\rm (i)}
Let the assumption be as in Lemma 8.1(i). Then, for
$\nu\leq r\leq d$, $\GrFB^{\nu}\CH^{r}(X;\QQ)$ contains
an uncountable number of classes with trivial de Rham invariant. 
\medbreak\noindent
{\rm (ii)}
Let the assumption be as in Lemma 8.1(ii). Then, for
$\nu\leq r\leq d$, $\GrFB^{\nu}\CH^{r}(X;\QQ)$ contains
an uncountable number of classes with trivial Mumford-Griffiths invariant. 
\endproclaim

\proclaim{Corollary 8.3} Let $X_o\in \C_{\ol{\QQ}}$ of dimension $d$,
and write $X=X_o\times_{\ol{\QQ}}\CC$. 
Assume given smooth curves $\Gamma_{1}/{\ol{\QQ}},
\ldots,\Gamma_{d}/{\ol{\QQ}}$ and a dominating morphism
$$
\Gamma_{1}\times \cdots \times \Gamma_{d} \to X_o.
$$
Further assume that $C(X)$ holds.
If $H^{0}(X,\Omega^\nu_{X/{\ol{\QQ}}})\ne 0$ for some $\nu\geq 2$,
then for $\nu\leq r\leq d$, \ 
$\GrFB^{\nu}\CH^{r}(X;\QQ)$ contains an uncountable number of classes 
with trivial Mumford-Griffiths invariant.
\endproclaim

Since any Abelian variety is dominated by a product of 
curves, we arrive at the following:

\proclaim{Corollary 8.3} Let $X_o$ be an Abelian variety over $\ol{\QQ}$, 
of dimension $d$ and write $X=X_o\times_{\ol{\QQ}}\CC$. . 
Then for any $2\leq \nu\leq r\leq d$, \
$\GrFB^{\nu}\CH^{r}(X;\QQ)$ contains
an uncountable number of classes with trivial Mumford-Griffiths invariant. 
\endproclaim

One can also arrive at similar results for Fermat hypersurfaces, using the 
fact that they are dominated by products of Fermat curves by 
\cite{KS}.

\newpage

\head Appendix: A conjectural overview\endhead

We would like to convey our thoughts
about the possibility of
a map from the space of algebraic cocycles 
$\Lambda_{\text{\rm alg}}^{r,\nu}(X/S)$ to
the space of Mumford-Griffiths invariants 
$\nabla J^{r,\nu}(X/S)\otimes \CC$, where
$\rho : X\to S$ is a smooth proper morphism
of smooth quasiprojective varieties over $\ol{\QQ}$.
In a nutshell, we believe that it ought to exist. We base
our observations on the situation with intersection 
cohomology considerations, and the various conjectures 
in the literature that seem to imply this.
\bigskip
\noindent
(I) The issue is whether there is a map from 
$$
\Lambda^{r,\nu}(X/S) := 
\hom_{\MHS}(\QQ(0),H^{\nu}(S,R^{2r-\nu}\rho_{\ast}\QQ(r))
\subset H^{\nu}(S,R^{2r-\nu}\rho_{\ast}\CC)
$$
at least on the space of algebraic
cocycles $\Lambda_{\text{\rm alg}}^{r,\nu}(X/S)$, (which
is Hodge conjecturally all of $\Lambda^{r,\nu}(X/S)$)
to the space of Mumford-Griffiths invariants $\nabla J^{r,\nu}(X/S)\otimes  \CC$. 
As already explained in \S8, if there is such a map, 
and if $S$ is affine, then such a map
must be injective. This is because, for $S$ affine, the de Rham
invariants $\nabla DR^{r,\nu}(X/S)\otimes\CC$ compute 
$H^{\nu}(S,R^{2r-\nu}\rho_{\ast}\CC)$, and the inclusion of
$\Lambda^{r,\nu}(X/S)$ in the
de Rham invariants factors through the map to $\nabla J^{r,\nu}
(X/S)\otimes \CC$. 
\bigskip
Whether or not there exists a map to the Mumford-Griffiths invariants,
we can still say that for $S$ affine, there is always a map
$$
\nabla J^{r,\nu}_{\text{\rm alg}}(X/S)\to 
\Lambda^{r,\nu}_{\text{\rm alg}}(X/S) \subset 
\Lambda^{r,\nu}(X/S),
$$
where 
$\nabla J^{r,\nu}_{\text{\rm alg}}(X/S)$ is the image
of $F^{\nu}\CH^{r}(X/S;\QQ)$ in $\nabla J^{r,\nu}(X/S)$.
One simply deduces this from the commutative diagram below,
where $\Phi^{r,\nu}$ is injective.
$$
\matrix
F^\nu \CH^r(X/S;\QQ) &@>{\lambda^{r,\nu}_{X/S}}>>& 
\Lambda^{r,\nu}(X/S) \\
&\\
\downarrow&&\downarrow\rlap{$\Phi^{r,\nu}$}\\
&\\
\nabla J^{r,\nu}(X/S)&\to& 
\nabla DR^{r,\nu}(X/S)\otimes\CC \\
\endmatrix
$$
The question then is, in more precise terms,
is the following.
\proclaim{Question A.0} For $S$ quasiprojective, does there
exist a map going the other way, viz.,
$$
\Lambda^{r,\nu}_{\text{\rm alg}}(X/S) \to \nabla J^{r,\nu}
(X/U)\otimes \CC,
$$
where $U/\ol{\QQ} \subset S$ is any
affine open subvariety? [Note:
If $S=U$ is affine, then this map must necessarily be injective,
if it exists.]
\endproclaim

\noindent
(II) Enter intersection cohomology: We work with this diagram:
$$
\matrix X&\hookrightarrow&\ol{X}\\
&\\
\rho\downarrow\quad&&\quad\downarrow\ol{\rho}\\
&\\
S&{\buildrel j\over\hookrightarrow}&\ol{S}
\endmatrix
$$
where $Y := \ol{X}\bs X$, $\Sigma := \ol{S}\bs S$ are NCD's, $\rho$ is smooth
and proper, and $\ol{X}$, $\ol{S}$ are smooth projective
over $\ol{\QQ}$.  We recall a
conjecture of Zucker and Brylinski.  First, let ${\Bbb V} = 
R^{m}\rho_{\ast}\CC$ over $S$. Put 
$$
{\Cal G}^{i} = \{g\in\Omega_{\ol{S}}^{i}\langle\log\Sigma\rangle\otimes 
\ol{\Cal V}\ |\ \text{\rm Res}_{\tilde{\Sigma}_{J}} g \in N_{J}V\},
$$
%Note $V$ is the fiber of ${\Cal V}$ over the point in $\Sigma$
%in question. It is a MHS.
where $\ol{\Cal V}$ is the privileged extension of
the bundle associated to ${\Bbb V}$ over $S$ to $\ol{S}$
(see \cite{B-Z} p. 63 for the definition of $\ol{\Cal V}$),
and $N_{J}$ for a multiindex $J$ corresponds to monodromy.
(See (III) below for a definition of the MHS $V$.)
All terminology in the conjecture below can
be found in \cite{B-Z}.
Then according to \cite{B-Z}, ``We suspect
that the means are available to verify:

\proclaim{Conjecture A.1} (\cite{B-Z}(3.19)) 
The complex ${\Cal G}^{\bullet}$
together with the [Hodge] filtration induced from (3-8) [of \cite{B-Z}], completes
[the middle perversity cochain complex] 
$\ul{IC}^{\bullet}(\ol{S},\Bbb {V}_{A})$ to a cohomological $A$-Hodge
complex of weight $m$, such that the Hodge structure it gives on
$IH^{i}(\ol{S},{\Bbb V})$ coincides with that
of (3.16) [of \cite{B-Z}].''
\endproclaim
For the benefit of the reader, we will clarify the statement
of Conjecture A.1. First of all, we make use of
the identification from \cite{C-K-S2},
$IH^{i}(\ol{S},{\Bbb V})  \simeq H^{i}_{(2)}(S,{\Bbb V})$,
where the former is intersection cohomology and
the latter is $L^{2}$-cohomology, which has a pure Hodge 
structure by \cite{C-K-S2}. This defines the Hodge structure
on $IH^{i}(\ol{S},{\Bbb V})$. (Alternatively, one has 
a Hodge structure on $IH^{i}(\ol{S},{\Bbb V})$ from the work
of M. Saito (\cite{MSa1}.) The privileged extension
$\ol{\Cal V}$ has corresponding Hodge filtration
subbundles $\ol{\Cal F}^{\bullet}\subset \ol{\Cal V}$,
which together with the monodromy weight filtration,
determine a MHS on the central fiber $V$. One has a corresponding
filtered subcomplex ${\Cal F}^{p}{\Cal G}^{i}  \subset
{\Cal G}^{i}$, where we observe that 
${\Cal F}^{p}{\Cal G}^{i}  \subset \Omega^{i}_{\ol{S}}\langle
\log \Sigma\rangle\otimes \ol{\Cal F}^{p-i}$, where the latter sheaf
is coherent. The import of Conjecture A.1 is that
$$
F^{\bullet}IH^{\bullet}(\ol{S},{\Bbb V}) \simeq \HH^{\bullet}
({\Cal F}^{\bullet}{\Cal G}^{\bullet}).
$$
Note that 
$$
{\Cal F}^{p}{\Cal G}^{i}\biggl|_{S} = \Omega_{S}^{i}\otimes
F^{p-i}R^{m}\rho_{\ast}\CC.
$$
Thus an immediate consequence of this conjecture 
(with $m = 2r-\nu$ and $p=r-\nu$) is the existence
of a map:
$$
\{H_{(2)}^{\nu}(S,R^{2r-\nu}\rho_{\ast}\CC)\}^{r,r} \to \nabla
J^{r,\nu}(X/U)\otimes\CC,\tag{A.2}
$$
where again $U/\ol{\QQ} \subset S$ is any affine open
subvariety.
%Regarding the map in (A.2): One has 
%$$
%\HH^{\nu}(F^{r}{\Cal S}^{\bullet}) \to 
%\HH^{\nu}(F^{r}[\Omega^{\bullet}_{\ol{S}}
%\langle\log \Sigma\rangle)\otimes \ol{\Cal V}]) \ {\buildrel
%\text{\rm GAGA}\over =}\  \HH_{\text{\rm Zar}}^{\nu}
%(\{F^{r}[\Omega^{\bullet}_{\ol{S}}
%\langle\log \Sigma\rangle)\otimes \ol{\Cal V}\}_{\text{\rm alg}}])
%\to \HH_{\text{\rm Zar}}^{\nu}(U,\{F^{r}[\Omega_{U}^{\bullet}\otimes 
%{\Cal V}\}_{\text{\rm alg}}]).
%$$
\bigskip
\noindent
\underbar{Remark A.2.1}. It would be nice if one could arrive
at and prove an analogous conjecture for 
$H^{\nu}(S,R^{2r-\nu}\rho_{\ast}\CC)$.
\bigskip

(III) Relating this to the original problem: We
will now assume given our map in (A.2) above.
Note that there is a natural map 
$$
H_{(2)}^{\nu}(S,R^{2r-\nu}\rho_{\ast}\CC) \to 
H^{\nu}(S,R^{2r-\nu}\rho_{\ast}\CC).
$$
Of course, in light of \cite{A}, one conjectures this to be a morphism
of mixed Hodge structures. Using
norm estimates on the weight filtration that arise
from the work of \cite{C-K-S1}, one argues
that:
\proclaim{Proposition A.3} There is a natural map
$$
H^{\nu}(\ol{S},j_{\ast}R^{2r-\nu}\rho_{\ast}\CC)
\to H_{(2)}^{\nu}(S,R^{2r-\nu}\rho_{\ast}\CC).
$$
\endproclaim

\demo{Proof} It suffices to construct a map:
$$
H^{\nu}(\ol{S},j_{\ast}R^{2r-\nu}\rho_{\ast}\CC)
\to IH^{\nu}(\ol{S},R^{2r-\nu}\rho_{\ast}\CC).
$$
Let
$\ul{IC}^{\bullet}(\ol{S},{\Bbb V})$ be the middle perversity 
intersection cochain complex, where ${\Bbb V} = R^{2r-\nu}\rho_{\ast}\CC$.
%This is well known
%to be a complex of fine sheaves. 
Then
$$
IH^{\nu}(\ol{S},{\Bbb V}) := {\Bbb H}^{\nu}\big(
\ul{IC}^{\bullet}(\ol{S},{\Bbb V})\big);
$$
The ``first'' hypercohomology spectral sequence (\cite{G-H})
computing ${\Bbb H}^{\nu}\big(
\ul{IC}^{\bullet}(\ol{S},{\Bbb V})\big)$ contains the $E_{2}$-term:
$$
{}^{\p}E_{2}^{\nu,0} := H^{\nu}\big(\ol{S},{\Cal H}^{0}(
\ul{IC}^{\bullet}(\ol{S},{\Bbb V})\big).
$$
Therefore, if there is a natural map
$$\tau:j_{\ast}{\Bbb V} \to {\Cal H}^{0}(\ul{IC}^{\bullet}(\ol{S},{\Bbb V})),$$
this, together with partial degeneration, leads to maps:
$$
H^{\nu}(\ol{S},j_{\ast}{\Bbb V}) \to {}^{\p}E_{2}^{\nu,0}
\to {}^{\p}E_{\infty}^{\nu,0}\hookrightarrow IH^{\nu}(\ol{S},{\Bbb V}),
$$
and hence Proposition A.3.
To construct $\tau$, by passing to a finite
cover of $S$, one can assume that the local monodromy 
transformations about
the NCD $\Sigma \subset \ol{S}$ are unipotent. 
In a classical
neighbourhood $U_{an}\cap  S \simeq (\Delta^{\ast})^{m}
\times \Delta^{\dim S -m}$.
There is no loss of generality in deleting factors $\Delta$.
Thus we may assume the local description $U_{an}\cap  S 
\simeq(\Delta^{\ast})^{m}$. Let $T_{j}$, $j=1,\ldots,m$
be the local monodromy transformations, and put $N_{j}
= \log T_{j}$. Let ${\Cal V} = {\Cal O}_{S}({\Bbb V})$,
and $\ol{\Cal V}$ the canonical extension over
$\Delta^{m}$. The stalk of $\ol{\Cal V}$ over $0 :=(0,\ldots0)\in 
\Delta^{m}$ will be denoted by $V$. Put
$$
B^{p}(N_{1},\ldots,N_{m};V)  = \bigoplus_{1\leq j_{1}<\cdots
< j_{p}\leq m}N_{j_{1}}N_{j_{2}}\cdots N_{j_{p}}V.
$$
Define the differential 
$$
(-1)^{s-1}N_{j_{s}} : N_{j_{1}}\cdots \widehat{N_{j_{s}}}
\cdots N_{j_{p}}V \to N_{j_{1}}N_{j_{2}}\cdots N_{j_{p}}V,
$$
on the various summands of $B^{p-1}(N_{1},\ldots,N_{m};V)$. This
turns $B^{\bullet}(N_{1},\ldots,N_{m};V)$ into a complex.
It is well known (\cite{C-K-S2}) that
$$
H^{\ast}_{(2)}(U_{an},{\Bbb V})\simeq
H^{\ast}\big(B^{\bullet}(N_{1},\ldots,N_{m};V)\big) \simeq 
IH^{\ast}(U_{an},{\Bbb V})=
({\Cal H}^{\ast}(\ul{IC}^{\bullet}(\ol{S},{\Bbb V})))_{0}.
$$
The existence of $\tau$ follows by noticing that the invariant cycles include
$$\big(j_{\ast}{\Bbb V}\big)_{0}\subset H^{0}
\big(B^{\bullet}(N_{1},\ldots,N_{m};V)\big).$$\qed

%(For example, if $\ol{S}$ is a curve and if $\nu = 1$, then
%one easily sees that $B^{0}(V) = V, \ B^{1}(V) = NV,\
%B^{\bullet \geq 2}(V) = 0$, hence ${\Cal H}^{i}(\ul{IC}^{\bullet}
%(\ol{S},{\Bbb V}) = 0$ for $i\geq 1$, and ${\Cal H}^{0}(\ul{IC}^{\bullet}
%(\ol{S},{\Bbb V}) = j_{\ast}{\Bbb V}$. Thus $H^{1}(\ol{S},
%j_{\ast}{\Bbb V}) \simeq IH^{1}(\ol{S},{\Bbb V})$ in this case.)

\enddemo

Using the natural map $R^{2r-\nu}\ol{\rho}_{\ast}\CC \to
j_{\ast}R^{2r-\nu}\rho_{\ast}\CC$, together with Proposition A.3,
we deduce:
\proclaim{Corollary A.4} There is a natural map
$$
H^{\nu}(\ol{S},R^{2r-\nu}\ol{\rho}_{\ast}\CC) \to
H_{(2)}^{\nu}(S,R^{2r-\nu}\rho_{\ast}\CC) \simeq IH^{\nu}
(\ol{S},R^{2r-\nu}\rho_{\ast}\CC).
$$
\endproclaim

In light of the mixed Hodge structure results in \cite{A},
it is natural to expect the following.

\proclaim{Conjecture A.5} The map in A.4 is
a morphism of MHS.
\endproclaim

Notice that $H^{\nu}(\ol{S},R^{2r-\nu}\ol{\rho}_{\ast}\CC)$
is the $E_{2}$ term of the Leray spectral sequence
for $\ol{\rho}$. 

\proclaim{Corollary A.6} Under the assumption of Conjecture A.5,
there is an induced morphism of Hodge structures:
$$
Gr_{\Cal L}^{\nu}H^{2r}(\ol{X},\CC) \to
H_{(2)}^{\nu}(S,R^{2r-\nu}\rho_{\ast}\CC),
$$
where  $Gr_{\Cal L}^{\nu}H^{2r}(\ol{X},\CC)$ is the
graded Leray filtration.
\endproclaim

\demo{Proof}
To see this, observe that $H_{(2)}^{\nu}(S,R^{2r-\nu}\rho_{\ast}\CC)$
has a pure Hodge structure of weight $2r$, and that by \cite{A},
$H^{p}(\ol{S},R^{q}\ol{\rho}_{\ast}\CC)$ has a mixed Hodge
structure of weight $\leq p+q$; moreover the differentials
of the Leray spectral sequence are morphisms of MHS
(see \cite{A}). Granted we accept Conjecture A.5,
it follows that the composite
$$
H^{\nu-2}(\ol{S},R^{2r-\nu+1}\ol{\rho}_{\ast}\CC)
@>d_{2}>> H^{\nu}(\ol{S},R^{2r-\nu}\ol{\rho}_{\ast}\CC) \to 
H_{(2)}^{\nu}(S,R^{2r-\nu}\rho_{\ast}\CC),
$$
is zero. Applying this same reasoning to the $d_{r}$'s for $r\geq 3$ leads to 
Corollary A.6.\qed
\enddemo

\proclaim{Corollary A.7} Let us assume that Conjectures A.1 and A.5 
hold. Then for $U/\ol{\QQ} \subset S$ affine, there is a map
$$
\{Gr_{\Cal L}^{\nu}H^{2r}(\ol{X},\CC)\}^{r,r} \to \nabla J^{r,\nu}
(X/U)\otimes\CC.\tag{$\ast$}
$$
\endproclaim

It is natural then to ask whether one can replace $\ol{X}$
in ($\ast$) by $X$?  The answer appears to be {\it yes,} provided
that one restricts to algebraic cocyles. The reason is
this. Notice that by purity of Hodge structures, and polarization, 
$Gr_{\Cal L}^{\nu}H^{2r}(\ol{X},\CC) \subset H^{2r}(\ol{X},\CC)$
as Hodge structures.  Let's
assume the Hodge conjecture for $Y$. More specifically,
if we write the NCD $Y$ in the form $Y = \bigcup_{j}Y_{j}$,
where $Y_{j}$ is smooth, then we are assuming the Hodge conjecture
for each $Y_{j}$. Then
the kernel of $H^{2r}_{\text{\rm alg}}(\ol{X})
\to H^{2r}(X)$ involves algebraic cocycles supported
on $Y$. But any algebraic cycle supported on
$Y$ goes to zero in $\CH^{r}(X)$. We want to make this
more precise. We need the following.

\proclaim{Lemma A.8} 
Consider the composite of the natural maps
$$
\tau\;:\; F_{\Cal L}^{\nu}H^{2r}(\ol{X},\QQ)  \to 
F_{\Cal L}^{\nu}H^{2r}(X,\QQ)  \to Gr_{\Cal L}^{\nu}H^{2r}(X,\QQ) \simeq
H^{\nu}(S,R^{2r-\nu}\rho_{\ast}\QQ),
$$
where the first map is the restriction and the isomorphism comes from 
the degeneration of the Leray spectral sequence for $\rho:X\to S$. 
Note that it is a map of Hodge structures and hence induces the map
$$
Gr_{\Cal L}^{\nu}H^{2r}(\ol{X},\QQ) \to 
W_{2r}H^{\nu}(S,R^{2r-\nu}\rho_{\ast}\QQ) = Gr^{2r}_{W}
H^{\nu}(S,R^{2r-\nu}\rho_{\ast}\QQ).
$$
Then the last map is a surjective morphism of HS.
\endproclaim

\demo{Proof} First of all, by
\cite{A}, the weight of $H_{c}^{\nu}(S,R^{2r-\nu}\rho_{\ast}\QQ)$
is $\leq 2r$. Hence by duality, the lowest weight space
of $H^{\nu}(S,R^{2r-\nu}\rho_{\ast}\QQ)$ is $2r$. That the
above is a morphism of MHS is the import of \cite{A}. This,
together with cohomology with compact supports yields the
exact sequence of MHS:
%Iverson: Cohomology of Sheaves
$$
H^{\nu-1}(\Sigma,R^{2r-\nu}\ol{\rho}_{\ast}\QQ) \to H_{c}^{\nu}(S,R^{2r-\nu}
\rho_{\ast}\QQ) \to H^{\nu}(\ol{S},R^{2r-\nu}\ol{\rho}_{\ast}\QQ).
$$
Since $Y\to \Sigma$ is a morphism of projective varieties,
the weight of $H^{\nu-1}(\Sigma,R^{2r-\nu}\rho_{\ast}\QQ)$ is
$\leq 2r-1$. 
%In \cite{A}, varieties do not have to be irreducible.
Thus we have the inclusion
$$
Gr_{W}^{2r}H^{\nu}_{c}(S,R^{2r-\nu}\rho_{\ast}\QQ)
\hookrightarrow Gr_{W}^{2r}H^{\nu}(\ol{S},R^{2r-\nu}\ol{\rho}_{\ast}\QQ).
$$
By compatibility of Leray spectral sequences together
with a weight argument similar to the proof
of Corollary A.6, this induces
an injection:
$$
Gr_{W}^{2r}H^{\nu}_{c}(S,R^{2r-\nu}\rho_{\ast}\QQ)
\hookrightarrow Gr_{W}^{2r}Gr_{\Cal L}^{\nu}H^{2r}(\ol{X},\QQ)
= Gr_{\Cal L}^{\nu}H^{2r}(\ol{X},\QQ).
$$
The Lemma A.8 now follows from duality.\qed
%cup product on E_{\infty} induced from that on E_{2}
%Also we use: $H^{i}_{c}(Y) \times H^{2\dim Y -i}(Y) \to \QQ$
%is perfect. This is just the pairing involving
%the extreme perversities in intersection (c0)-homology 
%theory with supports, and in a self dual locally constant system
%$L$, where in our case $L =R^{\bullet}\rho_{\ast}\QQ$.
%(e.g. See Borel's book). Now use Deligne's degeneration at $E_{2}$
%for the respective Leray spectral sequences.
\enddemo

\proclaim{Corollary A.9} {\rm (i)} 
For all $\nu\geq 0$, the kernel of $\tau$ in Lemma A.8 lies in the image 
$$
\big\{H^{2r}_{Y}(\ol{X},\QQ) + 
F_{\Cal L}^{\nu+1}H^{2r}(\ol{X},\QQ)\big\} \to H^{2r}(\ol{X},\QQ).
$$

{\rm (ii)} For all $\nu\geq 0$, the natural map
$$
F^{\nu}_{\Cal L}H^{2r}(\ol{X},\QQ) \to 
W_{2r}F^{\nu}_{\Cal L}H^{2r}(X,\QQ),
$$
is a surjective morphism of HS.
\endproclaim

\demo{Proof} Part (i): 
%Note that $W_{2r}H^{2r}(X,\QQ) =$ Image$\big(H^{2r}(\ol{X},\QQ) 
%\to H^{2r}(X,\QQ)\big)$,
%and that $Gr_{\Cal L}^{\nu}H^{2r}(X,\QQ) = H^{\nu}(S,R^{2r-\nu}
%\rho_{\ast}\QQ)$, hence the existence of the map in A.9. Next,
If $\xi\in F^{\nu}_{\Cal L}H^{2r}(\ol{X},\QQ)$  restricts
to a class in $F^{\nu+1}_{\Cal L}H^{2r}(X,\QQ)$, then by a repeated 
application of Lemma A.8, one can find $\xi^{\p}
\in F^{\nu+1}_{\Cal L}H^{2r}(\ol{X},\QQ)$ such that
$\xi -\xi^{\p}$ restricts to zero in $H^{2r}(X,\QQ)$.
This implies part (i). The proof of part (ii) is again
a repeated application of Lemma A.8, and
will be left to the reader.\qed
\enddemo

We now arrive at the following:

\proclaim{Theorem A.10} Let us assume that Conjectures A.1 and A.5 
hold, and that the Hodge conjecture holds. Then Question A.0
is answered in the affirmative.
\endproclaim

\demo{Proof} We have a commutative diagram
$$
\matrix
\CH^r(\ol{X}) &\to& \CH^r(X) \\
&\\
\downarrow\rlap{$c^r_{\ol{X}}$}&& \downarrow\rlap{$c^r_{X}$}\\
&\\
H^{2r}(\ol{X},\QQ(r)) &\to& H^{2r}(X,\QQ(r))\\ 
&\\
\bigcup &&\bigcup \\
&\\
F_{\Cal L}^{\nu}H^{2r}(\ol{X},\QQ(r)) &\to& 
F_{\Cal L}^{\nu}H^{2r}(X,\QQ(r))\\ 
\endmatrix
$$
where $c^r_{\ol{X}}$ and $c^r_{X}$ are the cycle class maps.
One can show that 
$c^r_{X}(F^\nu_B\CH^r(X/S))\subset F_{\Cal L}^{\nu}H^{2r}(X,\QQ(r))$
by the same argument as \cite{Sa3}, Prop.(2-1)
and that $\lambda^{r,\nu}_{X/S}$ is identified with the composite map:
$$\psi:  F^\nu_B \CH^r(X/S) @>{c^r_{X}}>> F_{\Cal L}^{\nu}H^{2r}(X,\QQ(r))
\to H^{\nu}(S,R^{2r-\nu}\rho_{\ast}\QQ(r)).$$
Our aim is to construct a natural map
Image$(\psi) \to \nabla J^{r,\nu}(X/S)\otimes\CC$ 
by using the map $(*)$ in  Corollary A.7. 
Put
$$
\CH^r(\ol{X}) \supset \Xi := (c^r_{\ol{X}})^{-1} 
\big(F_{\Cal L}^\nu H^{2r}(\ol{X},\QQ(r))\cap H^{2r}(\ol{X},\CC)^{r,r}\big) 
 \cap (j^*)^{-1}(F^\nu_B \CH^r(X/S)), 
$$
where $j^* : \CH^{r}(\ol{X}) \to \CH^{r}(X)$ is the restriction.
We have the following commutative diagrams:
$$
\matrix
\Xi &@>{j^*}>>& F^\nu_B\CH^r(X/S) \\
&\\
\downarrow\rlap{$c^r_{\ol{X}}$}&& \downarrow\rlap{$\psi$}\\
&\\
F_{\Cal L}^\nu H^{2r}(\ol{X},\QQ(r))\cap H^{2r}(\ol{X},\CC)^{r,r}
 &@>{\tau}>>& H^{\nu}(S,R^{2r-\nu}\rho_{\ast}\QQ(r)),\\
\endmatrix
$$

\vskip.25in

$$
\matrix
\Xi &@>{j^*}>>& F^\nu_B\CH^r(X/S) \\
&\\
\downarrow\rlap{$c^r_{\ol{X}}$}&& \downarrow\rlap{$\phi^{r,\nu}_{X/S}$}\\
&\\
F_{\Cal L}^\nu H^{2r}(\ol{X},\QQ(r))\cap H^{2r}(\ol{X},\CC)^{r,r}
 && \nabla J^{r,\nu}(X/S)\\ 
&\\
\downarrow&& \downarrow\\
&\\
Gr_{\Cal L}^\nu H^{2r}(\ol{X},\CC)^{r,r} &@>{(*)}>>& 
\nabla J^{r,\nu}(X/S)\otimes\CC,\\ 
\endmatrix
$$
Assuming the Hodge conjecture for the irreducible components of $Y$,
Corollary A.9 (i) implies that the kernel of $\tau$ in the former diagram is 
contained in the sum of $F_{\Cal L}^{\nu+1} H^{2r}(\ol{X},\QQ(r))$ and 
the image under $c^r_{\ol{X}}$ of the subgroup of 
$\CH^r(\ol{X})$ of the cycles supported on $Y$, which goes to zero under
the restriction $j^*$. Hence, by a diagram chase and Hodge conjecture
for $c^r_{\ol{X}}$, to construct a natural map
Image$(\psi) \to \nabla J^{r,\nu}(X/S)\otimes\CC$,
it suffices to show the surjectivity of  
$\Xi @>{j^*}>> F^\nu_B\CH^r(X/S)$.

Let $\alpha \in F^\nu_B\CH^r(X/S)$, and consider its
image $[\alpha] := c^r_{X}(\alpha) \in W_{2r}F_{\Cal L}^{\nu}H^{2r}(X,\QQ(r))$.
By Corollary A.9 (ii), there exists 
$$\ol{[\alpha]}\in F_{\Cal L}^{\nu}H^{2r}(\ol{X},\QQ(r))
\cap H^{2r}(\ol{X},\CC)^{r,r},$$
which restricts to $[\alpha]$. Next $\alpha$ is the restriction of a class
$\ol{\alpha}\in \CH^{r}(\ol{X};\QQ)$. Put $[\ol{\alpha}]
= c^r_{\ol{X}}(\ol{\alpha})$. Thus
$$
[\ol{\alpha}] - \ol{[\alpha]} \in \ker \big(H^{2r}(\ol{X},\QQ(r))
\to H^{2r}(X,\QQ(r))\big) \cap H^{2r}(\ol{X},\CC)^{r,r},
$$
and hence by the Hodge conjecture, 
$[\ol{\alpha}] - \ol{[\alpha]} = c^r_{\ol{X}}(\beta)$,
where $\beta\in \CH^r(\ol{X})$ is an algebraic cycle supported on $Y$.
Notice that
$$
\big(\ol{\alpha} -\beta\big)\biggl|_{X} = \alpha\in F^{\nu}_{B}\CH^{r}(X/S),
$$
as any cycle supported on $Y$ goes to zero under
the restriction $j^* : \CH^{r}(\ol{X}) \to \CH^{r}(X)$.
Furthermore 
$$c^r_{\ol{X}}(\ol{\alpha} -\beta)= \ol{[\alpha]} \in 
F_{\Cal L}^{\nu}H^{2r}(\ol{X},\QQ(r))\cap H^{2r}(\ol{X},\CC)^{r,r},$$
which proves the desired assertion.\qed

%Namely, if $v_{\nu}\in \ker\big(F_{\Cal L}^\nu H^{2r}(\ol{X},\QQ)^{r,r}
%\to Gr_{\Cal L}^\nu H^{2r}(\ol{X},\CC)^{r,r}\big)$, then
%Corollary A.9 (i) $\Rightarrow v_{\nu} = v_{\nu +1} + [\beta]$,
%where $v_{\nu+1} \in F_{\Cal L}^{\nu+1} H^{2r}(\ol{X},\QQ)^{r,r}$
%and $\beta \in (c^r_{\ol{X}})^{-1} (F_{\Cal L}^\nu H^{2r}(\ol{X},\CC))\cap 
%j^{-1}(F^\nu_B \CH^r(X/S))$. BY commutivity of the diagram,
%$[\beta] \mapsto 0\in \nabla J^{r,\nu}(X/S)\otimes \CC$,
%$\beta$ is supported on $Y$, and hence $j(\beta) = 0$.
%Further $v_{\nu+1} \in F_{\Cal L}^{\nu+1} H^{2r}(\ol{X},\CC)^{r,r}
%\mapsto 0\in Gr_{\Cal L}^\nu H^{2r}(\ol{X},\CC)^{r,r}$, hence
%to $0\in \nabla J^{r,\nu}(X/S)\otimes\CC$.
\enddemo

Recall from \S8 that in the case where $S$ is affine,
the composite
$$
\Lambda_{\text{\rm alg}}^{r,\nu}(X/S) \to \nabla J^{r,\nu}(X/S)
\otimes \CC \to \nabla DR^{r,\nu}(X/S)\otimes \CC,
$$
is injective.\footnote{Warning: The natural map $\nabla J^{r,\nu}(X/S)
\to \nabla DR^{r,\nu}(X/S)$ is not injective, as one can clearly
see from the product case $X = X_o\times S$.} Thus
the somewhat surprising conclusion from all of this, and from a conjectural
point of view, is that with regard to space of algebraic cycles,
one loses no information by passing from Mumford-Griffiths invariants
to de Rham invariants. Thus in particular,
this would enable us to replace classes with trivial
de Rham invariant, by trivial Mumford-Griffiths invariant in 
sections 7 and 8.

\bigskip
Finally, in the case $\dim S=1$, Theorem A.10 holds without
any conjectural assumption. By \cite{Z}, one has morphisms
of MHS:
$$
Gr_{\Cal L}^{1}H^{2r}(\ol{X},\CC) = 
H^{1}(\ol{S},R^{2r-1}\ol{\rho}_{\ast}\CC) \simeq H^{1}(\ol{S},
j_{\ast}R^{2r-1}\rho_{\ast}\CC) \hookrightarrow 
H^{1}(S,R^{2r-1}\rho_{\ast}\CC).
$$
Note that $H^{2}(S,R^{2r-2}\rho_{\ast}\CC) = 0$ if
$S$ is affine.

\newpage

\Refs

\item{[A]} D. Arapura, The Leray filtration is motivic.
{\tt http://www.math.purdue.edu/}$^{\sim}${\tt dvb/}.

\item{[As1]} M. Asakura M.,  
Motives and algebraic de Rham cohomology,
in: The Arithmetic and Geometry of Algebraic Cycles,
Proceedings of the CRM Summer School, June 7-19, 1998, Banff, Alberta, Canada
(Editors: B. Gordon, J. Lewis, S. M\"uller-Stach, S. Saito and N. Yui),
CRM Proceedings and Lecture Notes, {\bf 24} (2000), 133-155,
AMS.

\item{[As2]} Arithmetic Hodge structure and nonvanishing of the cycle class
of $0$-cycles, $K$-theory, {\bf 27} (2002), 273-280.

\item{[Bei]} A. Beilinson, Notes on absolute Hodge cohomology;
Applications of Algebraic $K$-Theory to Algebraic Geometry
and Number Theory, Contemporary Mathematics {\bf 55} (I), American
Mathematical Society,
Providence (1986), 35-68.

\item{[Bl]} S. Bloch;
Lectures on Algebraic Cycles,
Duke Univ. Math. Series {\bf 4} (1980)

\item{[Bo]} C. Borcea, Deforming varieties of 
$k$-planes of projective complete intersections,
 Pacific Journal of Math. {\bf 143} (1990), 25-36,
Sijthoff and Noordhoff.

\item{[B-Z]} J.-L. Brylinski and S. Zucker, An overview of
recent advances in Hodge theory, In: {\it Complex Manifolds,}
Springer-Verlag, New York, 1997, pp. 39-142.

\item{[Ca]} J. Carlson, Extensions of mixed Hodge structures,
 Journ\'ees de geom\'etrie alg\'ebrique d'Angers (1979), 107-127,
Sijthoff and Noordhoff.

\item{[C-K-S]} E. Cattani, A. Kaplan, W. Schmid, 1. Degeneration
of Hodge structures, Ann. of Math. {\bf 123} (1986), 457-535.

2. $L^{2}$ and intersection cohomologies for
a polarizable variation of Hodge structure, Inventiones math.
{\bf 87}, 217-252 (1987).

\item{[De1]} P. Deligne, Th\'eorie de Hodge II,
 Publ. Math. IHES,  {\bf 40} (1972), 5-57.
 
\item{[De2]} P. Deligne, Equations Diff\'erentielles \`a Points
Singuliers R\'eguliers, Lect. Notes in Math {\bf 163},
Springer-Verlag, 1970.

\item{[De3]} P. Deligne, Th\'eorem de L\`efschetz et crit\`eres de
d\'eg\'en\'erescence de suites spectrales, Publ. Math. I.H.E.S.
{\bf 35} (1969), 107-126.
 
 \item{[E-M]} Philippe Elbaz-Vincent and Stefan M\"uller-Stach,
 Milnor $K$-theory of rings, higher Chow groups and
 applications, Invent. math. {\bf 148}, (2002), 177-206.
 
\item{[E-P]} H. Esnault and K. H. Paranjape, Remarks on
absolute de Rham and absolute Hodge cycles,
C. R. Acad. Sci. Paris, t. {\bf 319}, S\'erie I,
(1994), 67-72.
 
\item{[EZ]} F. El Zein, Complexe dualisant et applications 
\`a la class fondamentale d'un cycle,
Bull. Soc. Math. France, M\'emoire {\bf 58} (1978).

\item{[F-G]} G. Frey and M. Jarden, Approximation theory and the
rank of Abelian varieties over large algebraic fields, Proc. London
Math. Soc. {\bf 28}(3), (1974), 112-128.

\item{[Gr]} M. Green,  Lectures at the CRM Summer School on 
``The arithmetic and Geometry of Algebraic Cycles", 
June 7-19, 1998, Banff, Alberta, Canada,
(handwritten notes).

\item{[Gr-Gr]} M. Green and P. A. Griffiths, Hodge-theoretic 
invariants for algebraic cycles, International Mathematics Research 
Notices 2003, No. 9, 477-510.

\item{[G-H]} P. A. Griffiths and J. Harris, {\it Principles
of Algebraic Geometry,} Wiley, New York, 1978.

\item{[Gri]} P. A. Griffiths, On the periods of
certain rational integrals: I and II, Annals of
Math. {\bf 90}(3), (1969), 460-541.

\item{[Ja]} U. Jannsen, Mixed sheaves and filtrations on Chow groups,
in: Motives (Editors: U. Jannsen, S. Kleiman, J.-P. Serre), 
Proceedings of Symposia in Pure Math. {\bf 55}, 245-302, AMS.

\item{[K]} M. Kashiwara and T. Kawai, The Poincar\'e lemma for 
variations of polarized Hodge structure, Publ. RIMS, Kyoto Univ. {\bf 
23} (1987), 345-407.

\item{[KS]} T. Katsura and T. Shioda, 
On Fermat varieties, Tohoku Math. J. {\bf 31}, (1979), 97-115.

\item{[KO]} N. Katz N. and T. Oda,
On the differentiation of De Rham cohomology classes with respect 
to parameters, J. Math. Kyoto Univ., {\bf 8} (1968), 199-213.

\item{[Ke1]} M. Kerr, Higher Abel-Jacobi maps for $0$-cycles.
Preprint, 2003.

\item{[Ke2]} M. Kerr, Private communication.

%\item{[Ki]} F. Kirwan, {\it An Introduction to Intersection Homology 
%Theory,} Pitman Research
%Notes in Mathematics Series {\bf 187}, Longman Scientific \& Technical
%(1988).

\item{[K1]} S. Kleiman, Algebraic cycles and the Weil conjectures, Dix
expos\'es sur la cohomologie des sch\'emas, North-Holland,
Amsterdam; Masson, Paris, 1968, 359-386.

\item{[K2]} S. Kleiman, The standard conjectures,
in: Motives (Editors: U. Jannsen, S. Kleiman, J.-P. Serre), 
Proceedings of Symposia in Pure Math. {\bf 55}, 3-20, AMS.

\item{[Ki]} J. R. King, Log complexes of currents and functorial properties
of the Abel-Jacobi map, Duke math. J. Vol. {\bf 50}, No. 1, 
(1983), 1-53.

\item{[Le1]} J. D. Lewis, {\it A Survey of the Hodge Conjecture,
CRM Monograph Series (Amer. Math. Soc.)} {\bf 10}, 1999.

\item{[Le2]} J. D. Lewis, A filtration on the Chow groups of a complex
projective variety, Compositio Math {\bf 128}, (2001),
299-322.

\item{[Le3]} J. D. Lewis, Cylinder homomorphisms and Chow groups, 
Math. Nachr. {\bf 160}, (1993), 205-221.

%\item{[Mac]} R. Macpherson, Global questions in the topology
%of singular spaces, Proc. Int. Congress of Math., Aug. 16-24, 1983, 
%Warszawa,
%213-235.

\item{[Mu]} D.  Mumford,
Rational equivalence of $0$-cycles on surfaces,
J. Math. Kyoto Univ., {\bf 9} (1969), 195-204.

%\item{[N]} M. V. Nori, Algebraic cycles and Hodge theoretic connectivity,
%Invent. Math. {\bf 111} (1993), 349-373.

\item{[R1]} A. A. Roitman,
Rational equivalence of zero-cycles, Math. USSR Sbornik 
{\bf 18}(4) (1972), 571-589.

\item{[R2]} A. A. Roitman, On $\Gamma$-equivalence of zero-dimensional 
cycles, Math. USSR Sbornik {\bf 15}(4) (1971), 555-567.

\item{[MSa1]} M. Saito, Mixed Hodge modules,
Publ. RIMS. Kyoto Univ. {\bf 26} (1990), 221-333.

\item{[MSa2]} M. Saito, Refined cycle maps,
in: Algebraic Geometry 2000, Azumino, 
Advanced Studies in Pure Math. {\bf 36} (2002), 115-144,
Mathematical Society of Japan, Tokyo.

\item{[MSa3]} M. Saito, On the formalism of mixed sheaves,
RIMS preprint.

\item{[MSa4]} M. Saito, Direct image of logarithmic
complexes, math.AG/0506020.

\item{[RS]} A. Rosenschon and M. Saito, Cycle map for
strictly indecomposable cycles,  Amer. J. Math. {\bf 125} (2003), 773-790.

\item{[Sa1]} S. Saito, Motives and filtrations on Chow groups,
Invent. Math. {\bf 125} (1996), 149-196.

\item{[Sa2]} S. Saito, Motives and filtrations on Chow groups, II,
in: The Arithmetic and Geometry of Algebraic Cycles,
Proceedings of the CRM Summer School, June 7-19, 1998, Banff, Alberta, Canada
(Editors: B. Gordon, J. Lewis, S. M\"uller-Stach, S. Saito and N. Yui),
NATO Science Series {\bf 548} (2000), 321-346,
Kluwer Academic Publishers.

\item{[Sa3]} S. Saito, Higher normal functions and Griffiths groups,
J. of Algebraic Geometry, {\bf 11} (2002), 161-201. 

\item{[V]} C. Voisin, Variations de structure de Hodge 
et z\'ero-cycles sur les surfaces
g\'en\'erals, Math. Ann., {\bf 299} (1994), 77-103.

\item{[S]} W. Schmid,  Variation of Hodge structure: The
singularities of the period mapping, Inventiones Math. {\bf 22}
(1973), 211-319.

\item{[Sch1]} C. Schoen, On Hodge structures and non-representability
of Chow groups, Compositio Math. {\bf 88} (1993), 285-316.

\item{[Sch2]} C. Schoen, Varieties dominated by product varieties,
International Journal of Mathematics {\bf 7}(4) (1996), 541-571.

\item{[So]}  C. Soul\'e, Op\'erations en $K$-th\'eorie alg\'ebrique, 
Canad. J. Math., {\bf 37}, No. 3, (1985), 488-550.

\item{[Sr]} V. Srinivas, Zero cycles on singular varieties, 
in: The Arithmetic and Geometry of Algebraic Cycles,
Proceedings of the CRM Summer School, June 7-19, 1998, 
Banff, Alberta, Canada
(Editors: B. Gordon, J. Lewis, S. M\"uller-Stach, S. Saito and N. Yui),
NATO Science Series {\bf 548} (2000), 347-382,
Kluwer Academic Publishers.

\item{[Z]} S. Zucker, Hodge theory with degenerating coefficients:
$L_{2}$ cohomology in the Poincar\'e metric, Ann. Math. {\bf 109} 
(1979), 415-476.

\endRefs

\enddocument

\bye